\documentclass[reqno]{amsart}

\usepackage{graphicx}

\theoremstyle{theorem}
\newtheorem{theorem}{Theorem}[section]

\newtheorem{lemma}[theorem]{Lemma}
\newtheorem{prop}[theorem]{Proposition}

\theoremstyle{definition}
\newtheorem{definition}{Definition}[section]
\newtheorem{remark}[definition]{Remark}

\numberwithin{equation}{section}

\newcommand{\CC}{\mathbb{C}}
\newcommand{\RR}{\mathbb{R}}
\newcommand{\QQ}{\mathbb{Q}}
\newcommand{\ZZ}{\mathbb{Z}}

\newcommand{\CP}{\mathbb{P}}

\newcommand{\CM}{\mathcal{M}}

\newcommand{\CL}{\mathcal{L}}

\newcommand{\G}{\gamma}
\newcommand{\PD}[1]{\frac{\partial}{\partial {#1}}}
\newcommand{\D}{\delta}

\begin{document}
\title[Clifford torus]
{Holomorphic discs, spin structures and \protect\\
       Floer cohomology of the Clifford torus  }

\author{Cheol-Hyun Cho}

\address{Department of Mathematics, University of Wisconsin, Madison, WI
53706, Current address: Department of Mathematics, Northwestern University,
Evanston, IL 60208, cho@math.wisc.edu }

\begin{abstract}
We compute the Bott-Morse Floer cohomology of the Clifford torus
in $\CP^n$ with all possible spin-structures. Each spin structure
is known to determine an orientation of the moduli space of
holomorphic discs, and we analyze the change of orientation
according to the change of spin structure of the Clifford torus.
Also, we classify all holomorphic discs with boundary lying on the
Clifford torus by establishing a Maslov index formula for such
discs. As a result, we show that in odd dimensions there exist two
spin structures which give non-vanishing Floer cohomology of the
Clifford torus, and in even dimensions, there is only one such
spin structure. When the Floer cohomology is non-vanishing, it is
isomorphic to the singular cohomology of the torus (with a Novikov
ring as its coefficients).  As a corollary, we prove that any
Hamiltonian deformation of the Clifford torus intersects with it
at least at $2^n$ distinct intersection points, when the
intersection is transversal.

We also compute the Floer cohomology of the Clifford torus with
flat line bundles on it and verify the prediction made by Hori
using a mirror symmetry calculation.
\end{abstract}
\maketitle

\bigskip

\section{Introduction}

 The Floer homology of Lagrangian
intersection was first defined by Floer \cite{floer:msli88} and
since then, it is emerging as a powerful technique in symplectic
geometry. It has received much more attention after Konsevich
\cite{konsevich:hams95} proposed a homological Mirror symmetry
conjecture to use Floer homology in the context of $A_{\infty}$
category that Fukaya introduced\cite{fuk1}. The construction of Floer
homology has been generalized and applied to the problems in
symplectic geometry by Oh \cite{oh:fclipd93}, \cite{oh:fclipd293}
,\cite{oh:fsm96}, and recently, it was studied in complete
generality and its obstruction to the well-definedness of Floer
homology was established by Fukaya, Oh, Ohta and Ono in
\cite{fooo:lift00}. But computing actual Floer homologies is still
a difficult task, since one has to analyze all holomorphic strips
with boundary lying on two Lagrangian submanifolds. The
construction of Floer homology in the Bott-Morse setting is a big
step forward in this respect as in Morse theory.

In this paper, we compute Floer cohomology of the Clifford torus
$T^n$ in $\CP^n$ in the Bott-Morse setting. There are two main
issues in the computation. The first one is an orientation
problem. Floer and Oh defined Floer homology with $\ZZ /
2\ZZ$-coefficients. In \cite{fooo:lift00}, Fukaya, Oh, Ohta and
Ono developed a way to give an orientation of the moduli space of
holomorphic discs(and strips). This orientation depends on the
(relative) spin structure of a Lagrangian submanifold. Hence, one
can define Floer homology with $\ZZ$ or $\QQ$-coefficient. A {\em
spin structure} of an oriented vector bundle $E$ over $X$ can be
understood as a homotopy class of a trivialization of $E$
 over the 1-skeleton of $X$ which can be extended to
the 2-skeleton of $X$. It is already observed in
\cite{fooo:lift00} that different homotopy classes of the
trivialization of a certain bundle will reverse orientation of the
moduli space, with an example in the case of Maslov index 0 disc.
Here we give a proof of this observation in general by using the
Index theorem of Silva in \cite{silva:psfh97} (see Theorem
\ref{sgnchange}).

The Clifford torus is an interesting example since it has $2^n$
different spin structures. These $2^n$ spin structures can
possibly give rise to $2^n$ different Floer cohomologies. Or more
generally, we can consider Floer cohomology of $T^n$ with flat
line bundles on it. It may be considered as an advantage in that
we may exploit this freedom to define non-vanishing Floer
homology, if possible.

But for the Clifford torus, one can choose a natural spin
structure which we will call the {\em standard} spin structure.
Under the standard spin structure, it is not hard to determine the
orientation of the moduli space as described in \cite{fooo:lift00}
(see section \ref{sec:torusori}). For the other spin structures,
we will determine the orientation of the moduli space by studying
the change of orientation with respect to the change of spin
structures. Hence, we can determine orientations needed to define
Floer boundary operator for any spin structures of the Clifford
torus.

The second issue is to classify all the holomorphic discs with
boundary lying on a Lagrangian submanifold. For that purpose, we
prove the Maslov index formula (Theorem \ref{indexformula}) and
classify all the holomorphic discs with boundary on the Clifford
torus with any Maslov index. In this case, any non-constant
holomorphic disc has positive Maslov index, in which case
Bott-Morse Floer homology is rather easy to define.

By the classification theorem, Theorem \ref{classify}, we can
explicitly calculate Bott-Morse Floer boundary operators for Floer
cohomology. It turns out that among $2^n$ spin structures, for
$n=dim(L)$ even, there is only one spin structure which gives
non-vanishing Floer cohomology. And for $n$ odd, there exists two
spin structures which gives non-vanishing Floer cohomology. And
when it is non-vanishing, $HF(T^n,T^n;\Lambda_{nov})$ as a
$\Lambda_{nov}$-module is isomorphic to the singular cohomology
with $\Lambda_{nov}$-coefficient.

One immediate corollary of the latter result is that intersection
between the Clifford torus and its Hamiltonian deformation must
have at least $2^n$ distinct points when the intersection is
transversal. In particular, the Clifford torus must intersect
{\it any} Hamiltonian deformations thereof.  While we are in the 
preparation of
the thesis \cite{cho:hsf03}, we have learned from Oh that  this latter intersection
result was also proved by Biran-Entov-Polterovich \cite{bep} using
a completely different method without using the Floer homology. 

As an application to physics, one can compute a D-brane Floer
cohomology(Floer cohomololy with flat line bundle on the Lagrangian
submanifold). The homological mirror symmetry conjecture is about
Calabi-Yau manifolds, but, its extension to Fano case has been
studied by Hori \cite{hori:lmsd00}. With minor modification from
our calculation of Floer cohomology, we can compute compute
D-brane Floer cohomology. As a result, we found $(n+1)$ flat line
bundles with specific holonomies over the Clifford torus whose Floer cohomology is
non-vanishing, which has been predicted by Hori \cite{hori:lmsd00},
Hori-vafa \cite{hori:ms02}
by B-model calculation.

More generally, the prediction by K. Hori about the Floer
cohomology of Lagrangian torus fibers of Fano toric manifolds is
that the Floer cohomology of all the fibers vanish except at a
finite number of base points in the momentum polytope that are
critical points of the super-potential of the Landau-Ginzburg
mirror to the toric manifold. We generalize the scheme used in the
paper to this case and will prove the exact correspondence in
\cite{co}.

This is the simplified version of the author's Ph. D. thesis in
the University of Wisconsin-Madison.

We would like to thank Yong-Geun Oh for helpful discussions and
invaluable support.

\section{The Maslov index}\label{maslovindex}

In this section, we recall basic definitions including the Maslov
index of a map and its generalization in the case that the domain
of a map is a smooth Riemann surface with boundary.

Let $(M,\omega)$ be a 2$n$-dimensional compact symplectic
manifold. Let $L$ be a Lagrangian submanifold. There are two
homomorphisms $I_{\omega},\mu$ on $\pi_2(M,L)$ defined as follows.
The symplectic energy $I_{\omega} :\pi_2(M,L) \to \RR$ is defined
as $$I_{\omega} = \int_{D^2} w^* \omega.$$ The Maslov index $\mu
:\pi_2(M,L) \to \ZZ$ is defined as follows: We first consider the
Lagrangian Grassmannian $\Lambda(\CC^n)$ consisting of all
$n$-dimensional $\RR$ linear subspaces $V$ of $\CC^n$ such that
the standard symplectic form $\omega_0$ of $\CC^n$ vanishes on
$V$. The unitary group $U(n)$ acts transitively on
$\Lambda(\CC^n)$ and the isotropy group is $O(n)$. Therefore, we
have $\Lambda(\CC_n) \cong U(n) /O(n).$ Each Lagrangian plane can
be written as $A \cdot \RR^n$ for some $A \in U(n)$ and two such
matrices $A_1$,$A_2$ define the same plane if and only $A_1 \cdot
\overline{A}_1^{-1} = A_2\cdot \overline{A}_2^{-1}$. By
Proposition 4.2 \cite{oh:rhp95}, the map
$$B:\Lambda(\CC_n) \to \widetilde{\Lambda(\CC_n)}:\; A \mapsto A \cdot \overline{A}^{-1} = A\cdot A^t$$
is a diffeomorphism for
$$\widetilde{\Lambda(\CC_n)} = \{ D\in GL(n;\CC)| D\overline{D} = Id, D = D^t\}$$
Now, for any loop $\gamma : S^1 \to \Lambda(\CC_n)$, the {\bf
Maslov index} of a loop $\gamma$ is defined to be the degree of
the map $\phi = det \circ B \circ \gamma:S^1 \to U(1)$.

Now let $w:(D^2,\partial D^2) \to (M,L)$ be a smooth map
representing the homotopy class $\beta \in \pi_2(M,L)$. Then we
can find a unique trivialization (up to homotopy) of the pull-back
bundle $w^*(TM) \simeq D^2 \times \CC^n$ as a symplectic vector
bundle. The trivialization defines a map from $\gamma :\partial
D^2 \to \Lambda(\CC^n)$ and we define
$$\mu(w) :=\mu(\gamma) \in \ZZ.$$
It is independent of the trivialization. We will call $\mu(\beta)
= \mu(w)$ the {\bf Maslov index} of $\beta$. The {\em minimal
Maslov number} $\Sigma_L$ is the positive generator of the abelian
subgroup $[\mu|_{\pi_2(M,L)}] \subset \ZZ$.
\begin{definition}
A Lagrangian submanifold $L$ is said to be {\em monotone} if there
exists $c>0$ independent of $\beta \in \pi_2(M,L)$ such that
$$\mu(\beta) = c I_{\omega}(\beta).$$
\end{definition}
Let $\Sigma$ be a smooth Riemann surface with boundary. We will
denote by $R_0,\cdots,R_h$ the connected components of $\partial
\Sigma$, with orientation induced by the orientation of $\Sigma$.
We assume that the number of boundary components $h$ is nonzero.
Let $w:(\Sigma,\partial \Sigma) \to (M,L)$ be a smooth map with
$w(\partial \Sigma) \subset L$. Then we can also define the Maslov
index of the map $w$ as follows (see \cite{katz:egsm01}).

Let $E$ be the complex vector bundle $w^*TM$, and let $E_R$ be the
Lagrangian subbundle $w|_{\partial \Sigma}^*TL$. Since any complex
vector bundle over a Riemann surface with nonempty boundary is
trivial, we may fix the trivialization of $E$ as $\Phi:E \cong
\Sigma \times \CC^n$. Then, for each boundary component $R_i$, we
have a map $\gamma_i : S^1 \to \Lambda(\CC^n)$. Let
$\mu(\Phi,R_i)=\mu(\gamma_i)$. We define the {\bf Maslov index} of
the map $w$ as $$\mu(w) = \mu(\Phi,w) = \sum_{i=0}^{h}
\mu(\Phi,R_i)$$
\begin{prop}[Katz-Liu~\cite{katz:egsm01} Proposition 3.3.6]\label{maslovind}
The Maslov index defined above is independent of the choice of
trivialization $\Phi : E \cong \Sigma \times \CC^n$. $\Box$
\end{prop}

\section{The Clifford torus}
We follow the description of the Clifford torus given in
\cite{oh:fclipd93}. Consider the isometric embedding
 $$ T^{n+1} := \underbrace{ S^1(\frac{1}{\sqrt{(n+1)}}) \times \cdots \times
 S^1(\frac{1}{\sqrt{(n+1)}})}_{n+1 \,\textrm{times}}
 \hookrightarrow S^{2n+1}(1) \subset \CC^{n+1}$$
This embedding is Lagrangian in $\CC^{n+1}$, and the standard
action by $S^1$ on $\CC^{n+1}$ restricts to both the above torus
and $S^{2n+1}(1)$. By taking the quotients by this action, the
torus $T^n := T^{n+1}/S^1$
 in $\CP^n = S^{2n+1}(1)/S^1$ is Lagrangian submanifold. This torus is
 a minimal submanifold in Riemannian geometry; it is called the
{\bf Clifford torus} in $\CP^n$. For the case $n=1$, $T^1$ is
nothing but the great circle in $\CP^1$

\begin{prop}[Oh~\cite{oh:fclipd93} Proposition 2.4]\label{torusmonotone}
The above Clifford torus $T^n \subset \CP^n$ is monotone with
respect to the standard symplectic structure on $\CP^n$. $\Box$
\end{prop}
\begin{proof} We first describe the homotopy classes in
$\pi_2(\CP^n,T^n)$. We have the homotopy exact sequence,
$$\to \pi_2(T^n) \to \pi_2(\CP^n) \stackrel{i}{\to} \pi_2(\CP^n,T^n)
\stackrel{j}{\to} \pi_1 (T^n) \to \pi_1(\CP^n) \to $$
 with $\pi_2(T^n) \cong 0$ and $\pi_1(\CP^n) \cong 0$.
We have
$$\pi_2(\CP^n,T^n) \cong \pi_2(\CP^n)\oplus \pi_1 (T^n)$$
since the boundary map has an obvious right inverse.

For $0 \leq i \leq n$, let $b_i$ be the holomorphic disc
$$b_i=[\underbrace{1:\cdots:1}_{i}:z:1:\cdots:1]$$
We will denote their homotopy classes as $\beta_i = [b_i]\in
\pi_2(\CP^n,T^n)$. These are discs with the Maslov index 2, and we
will call them standard discs. Later, we will show that any 
holomorphic disc of Maslov index 2 with boundary lying on $T^n$
is in fact one of the standard discs up to an automorphism of a disc.

Now we want to show that the spherical homotopy class
$i(\pi_2(\CP^n)) \subset \pi_2(\CP^n,T^n)$ can be obtained as a
sum of $b_i$. For the generator $\alpha \in \pi_2(\CP^n)$ with
$c_1(\alpha) = n+1$ where $c_1$ is the first chern class of the
tangent bundle of $\CP^n$, it is known that the Maslov index of
$i(\alpha)$ is actually $2c_1(\alpha)=2(n+1)$. Now, note that
$$j(\beta_0 + \beta_1 + \cdots + \beta_n) =0 \in \pi_1 (T^n)$$ and
$\mu(\beta_0 +  \cdots + \beta_n) = 2n+2$. Since the homotopy
sequence is exact, $(\beta_0 +  \cdots + \beta_n)$ lies in the
image of the map $i:\pi_2(\CP^n) \to \pi_2(\CP^n,T^n)$. Hence, we
have $$ i(\alpha) = \beta_0 + \beta_1 + \cdots + \beta_n.$$

Then, it is easy to show that the Lagrangian submanifold $T^n$ is
monotone: If $I_{\omega}(\beta_i) = c \mu(\beta_i)$ for all $i$
for fixed $c$, then
$$I_{\omega}(i(\alpha)) = (n+1)I_{\omega}(\beta_i) = (n+1) c \mu(\beta_i)
=c \mu(i(\alpha)).$$ This proves that the
Clifford torus is monotone.
\end{proof}

\section{Bott-Morse Floer cohomology}\label{bott}
We review the construction of Bott-Morse Floer cohomology
$HF(L,L)$ following \cite{fooo:lift00}.
 There is a canonical isomorphism $HF(L,L) \to HF(L,\phi(L))$.
\begin{definition}[Novikov ring]\label{novikov}
We consider the formal (countable) sum $\sum_{i=0}^{\infty} c_i
e^{d_i}$ such that
$$ c_i \in \QQ,\,\, d_i \in \ZZ,\,\, lim_{i \to \infty} d_i = \infty$$
The totality of such formal sums becomes a ring, and we denote
this ring by $\Lambda_{nov}$. We consider $\sum_{i} c_i e^{d_i}$
with $d_i \geq 0$ in addition and denote it by $\Lambda_{0,nov}$.
Here we set the degree of $e$ to be $2$.
\end{definition}
\begin{remark}
Since we only consider monotone Lagrangian submanifolds, we
do not need to include the energy term here.
\end{remark}

To construct Floer cohomology in this case, we need a cochain
complex which represents cohomology theory of L. For a given
(n-k)-dimensional geometric chain $[P,f]$, we consider the current
$T([P,f])$ which is defined as follows: The current $T([P,f])$ is
an element in $D'^k(M;\RR)$ where $D'^k(M;\RR)$ is the set of
distribution valued k-forms on $M$ : For any smooth (n-k)-form
$\omega$, we put
\begin{equation}\label{current}
\int_M T([P,f]) \wedge \omega = \int_P f^* \omega
\end{equation}
This defines a homomorphism $$ T: S_{n-k}(M;\QQ) \to D'^k(M;\RR)$$
where $S_{n-k}(M;\QQ)$ is the set of all (n-k) dimensional
geometric chains with $\QQ$-coefficient. Let
$\overline{S}^k(M,\QQ)$ be the image of the homomorphism $T$. Then
we take a countably generated subcomplex $C(L;\QQ)$ of
$\overline{S}^k(M,\QQ)$ such that the cohomology of $C$ is
isomorphic to the cohomology of $H^*(M,\QQ)$. Since we consider
the elements in the image of $T$, if the image of the map $f$ of
the geometric chain $[P,f]$ is smaller than expected dimension,
then it gives 0 as a current. This fact will be used crucially
later on.


We recall the definition of the compactified moduli space of holomorphic discs
(See \cite{fooo:lift00} for details).
\begin{definition}[\cite{fooo:lift00}]
Let $\beta \in \pi_2(M,L)$ and denote by $\CM_m(\beta)$ the set of
all isomorphism classes of genus 0 stable maps from open curve with
(m,0) marked points $((\Sigma,\vec{z}),w)$ such that
$w_{*}([\Sigma])=\beta $.
Also denote by $\CM_m^{reg}(\beta)$ the subset of $\CM_m(\beta)$
with $\Sigma = D^2$. 
\end{definition}

For the analysis of orientations, we define the moduli space of 
holomorphic discs without compactification.
\begin{definition}\label{modulinoquo}
For a given homotopy class $\beta \in \pi_2(M,L;x)$,
We define
$$ \widetilde{\CM} (L,J : \beta) = \{ w : D^2 \to M |
\overline{\partial}_J w = 0, w ( \partial D^2) \subset L, [w] = \beta \}.$$
\end{definition}
We also similarly define $\widetilde{\CM}_n(L,J : \beta)$ to be 
the moduli space of holomorphic discs with $n$ marked points.
We will abbreviate $ \widetilde{\CM}_n (L,J : \beta) $ as 
$ \widetilde{\CM}_n (\beta)$ from now on.
The group $PSL(2:\RR)=Aut(D^2,j_D)$ acts on $\widetilde{\CM} (\beta)$ by
$\phi \cdot w = w \circ \phi^{-1}$ for $\phi \in PSL(2:R)$ and
$w \in \widetilde{\CM} (\beta)$ and it acts on a marked point $z_i$
as $\phi(z_i)$.
Then, we have $$\widetilde{\CM}_n (\beta)/PSL(2:\RR) \cong \CM_n^{reg}(\beta)$$

Now we recall the definition of Floer coboundary operator.
\begin{definition}
For a geometric chain $[P,f]\in C^*(L:\QQ)$, define
\begin{equation}
\begin{cases}
\D_{\beta}([P,f])=(\CM_2(\beta) \,_{ev_1}\times_f P, ev_0)
 \,\,\,\,\textrm{for }\, \beta \neq 0,\\
\D_0([P,f])=(-1)^n[\partial P,f]
\end{cases}
\end{equation}
\end{definition}
\begin{remark}
Well-definedness of this fiber product is rather technical,
because the moduli space possibly has codimension 1 corners and
the product is defined in the chain level. One need
Smooth-correspondence developed in \cite{fooo:lift00} appendix
$A$. However, in our later calculation, only non-trivial fiber
product occurs for the moduli space of Maslov index 2 discs, in
which case, the moduli space is closed (without boundary) since
the homotopy class is minimal. And we use the spectral sequence to
compute the Floer cohomology, therefore, after the first step, we
can work on the homology level.
\end{remark}


\begin{theorem}[\cite{fooo:lift00} Proposition 13.16]\label{dim}
For $[P,f]\in C^k(L:\QQ)$ ,
$$\delta_\beta([P,f]) \in  C^{k-\mu(\beta)+1}(L:\QQ).$$
\end{theorem}

Now we define our coboundary operator $\D$ on
$C^*(L;\Lambda_{nov})$ by extending the following boundary
operater linearly over $\Lambda_{nov}$
$$\delta([P,f])=\sum_{\beta \in \pi_2(M,L)}
\delta_\beta([P,f]) \otimes e^{\frac{\mu(\beta)}{2}}$$
\DeclareGraphicsExtensions{.eps,.jpg}
\begin{theorem}\label{squarezero}
$$\delta \circ \delta =0.$$
\end{theorem}
\begin{remark}
This is a combination of arguments used in \cite{fooo:lift00} Theorem 6.24
and [addenda,O1]. The proof in [FOOO] deals with the case when
$L$ is un-obstructed, while in this case obstructuion cycle is 
a multiple of the fundermental class $[L]$.
\end{remark}
\begin{proof}
It is enough to show that $\delta \circ
\delta ([P,f])=0$. Now,
\begin{equation}\label{dd0}
\delta \circ \delta ([P,f]) = \sum_{A \in \pi_2(M,L)} \sum_{A_1 +
A_2 = A} \delta_{A_1} \circ \delta_{A_2} [P,f] \otimes
e^{\frac{\mu(A)}{2}}
\end{equation}
We consider the geometric chain $(\CM_2(A) \,_{ev_1}\times_f P,
ev_0)$. Note that we consider not the moduli space itself but its
image under evaluation map.

\begin{figure}\begin{center}
\includegraphics[height=3in]{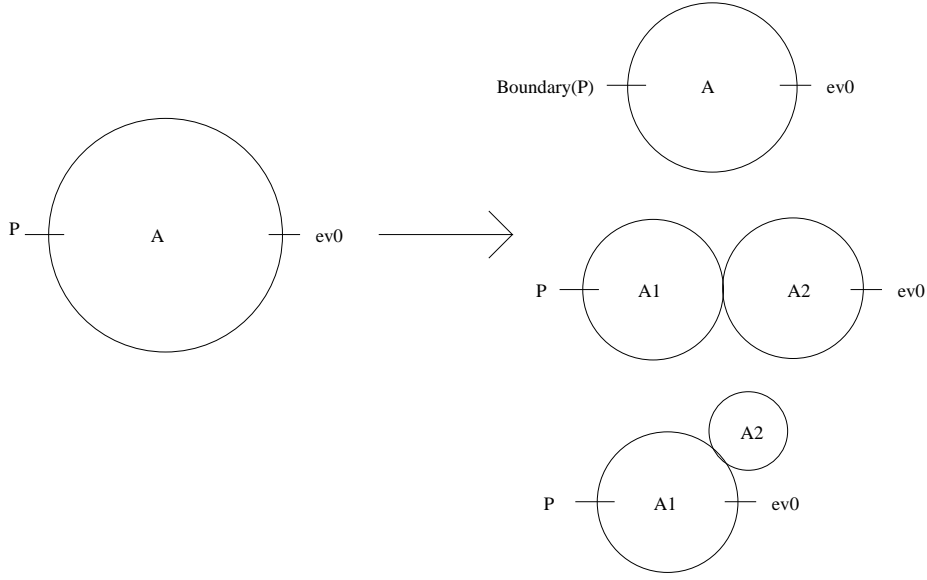}
\caption{Boundaries of $ev_0(\CM_2(A) \,_{ev_1}\times_f P)$}
\label{fig1}
\end{center}
\end{figure}

Now, we can describe the boundary components of the image of the
chain as in Fig \ref{fig1}. First there is a component corresponding to boundary of the
chain $P$, and a splitting of the moduli space of $A$,  and also there
is a component with a disc or sphere bubble. A component with
sphere bubble has at least codimension 2, hence it causes no
trouble. But a component with a disc bubble has codimension 1. But
in this Morse-Bott setting, in the case of monotone Lagrangian
submanifolds, image of such a component with a disc bubble does
not appear as a codimension 1 boundary as follows: (Generally, disc bubbling
phenomenon causes trouble defining $HF(L_0,L_1)$ for two different
Lagrangian submanfolds $L_0,L_1$. But in the case $L_1 =
\phi(L_0)$, or $L_0 = L_1$ disc bubbling with positive maslov index does not
cause much trouble defining Floer cohomology. See Proposition 7.3
\cite{fooo:lift00}). Basically, we only consider the image under
the evaluation map and we claim that the image is always of
codimension 2 or more.

As in Fig \ref{fig1}, if the disc $A$ splits, we call the component meeting the chain
$P$ as $A_1$ and the other component as $A_2$.

First, consider the case that $\mu(A_1)\neq 0$. Note that
$\mu(A_2) \geq 2$ since homotopy class of a bubble is always
non-trivial and the Lagrangian submanifold is orientable. Then the
image under the evaluation map of such a component is contained in
$(\CM_2(A_1) {\,}_{ev_1}\times_f P, ev_0)$ whose chain dimension is
$(n-k) + \mu(A_1) -1$. But the original chain $(\CM_2(A)
{\,}_{ev_1}\times_f P, ev_0)$ has chain dimension $(n-k) + \mu(A)
-1$. Since $\mu(A) \geq \mu(A_1) +2$, the image is of codimension 2 
or more as claimed.
\begin{figure}
\begin{center}
\includegraphics[height=2in]{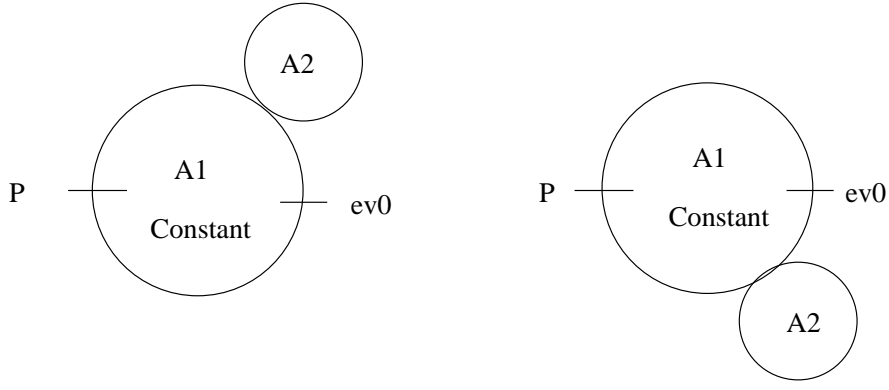}
\caption{Cancellation of disc bubbles when $\mu(A_1)=0$}
\label{fig2}
\end{center}
\end{figure}

Now, we consider the case that $\mu(A_1) = 0$. Then, actually
there should be pairs of bubbles occur as in the Fig \ref{fig2}.
Since $A_1$ is a constant holomorphic disc, its image under the
evaluation maps are the same. But these two components have different
orientation because of the ordering of the marked points.
Therefore these two bubbles cancel out each other. This proves
that the image of the component with a disc bubble will not give a
codimension 1 boundary.

The remaining boundary components can be written as follows.
\begin{eqnarray*}
\delta_0 \circ \delta_{A} [P, f] &=& \delta_0
(\CM_2(A) \,_{ev_1}\times_f P) \\
&=&(-1)^n \, \partial ( \CM_2(A) \,_{ev_1}\times_f P) \\
&=&(-1)^n (\partial \CM_2(A) \,_{ev_1}\times_f P) \bigsqcup (-1)^n
(-1)^{dim \CM_2(A) + n}
(\CM_2(A) \,_{ev_1}\times_f \partial P ) \\
&=& (-1)^n (\partial \CM_2(A) \,_{ev_1}\times_f P) \bigsqcup
(-1)^{n+n+\mu (A) +2 -3+n +n}
(\CM_2(A) \,_{ev_1}\times_f \delta_0 P ) \\
&=& (-1)^{n +n+1} ((\CM_2(A_1) \,_{ev_1} \times_{ev_0} \CM_2(A_2)
\,_{ev_1}\times_f P ) \bigsqcup
(-1) (\CM_2(A) \,_{ev_1}\times_f \delta_0 P ) \\
&=& (-1) \delta_{A_1} \circ \delta_{A_2} [P,f] \bigsqcup
(-1) (\CM_2(A) \,_{ev_1}\times_f \delta_0 P ) \\
\end{eqnarray*}
The third equality is from the formula~\ref{eq:ori1} in
section~\ref{sec:ori} and the fifth equality is from the
formula~\ref{eq:ori2} in the same section. Hence it proves that
$$\delta \circ \delta ([P,f])= \delta_0 \circ \delta_{A} [P, f] +
\delta_{A_1} \circ \delta_{A_2} [P,f] + \delta_{A} \circ \delta_0
[P, f]=0.$$
\end{proof}


We also recall the construction of the spectral sequence which
converges to $HF(L,L)$ for the monotone Lagrangian submanifold
$L$. Existence of the spectral sequence was first observed by Oh
\cite{oh:fsm96}. When the Lagrangian submanifold is monotone, we
have the minimal Maslov index $\Sigma_L$ of $L$. Let $\delta_i$ be
the formal sum of $\delta_{\beta}$ with the Maslov index $i$.
Then, we have
$$\delta = \delta_0 + \delta_{\Sigma_L}\otimes e^{\Sigma_L/2} + \delta_{2\Sigma_L}
\otimes e^{\Sigma_L} + \cdots $$ This filtration basically gives
the spectral sequence of the Floer cohomology. The spectral
sequence will start from cochain complex and the boundary maps in
$E_i^{*,*}$ will be $\delta_{2i-2}$. For $L$ monotone, filtration
by energy and that by Maslov index are equivalent. But in general,
one should consider filtration by energy (see \cite{fooo:lift00}).

\begin{theorem}[\cite{fooo:lift00} Theorem 6.13]\label{spectral}
There exists a spectral sequence with
$$E_2^{p,q} \cong \oplus_{i=q}^{[\frac{p}{2}]} (H^{p-2i}(L:\QQ)\otimes e^i) \cong
(H^*(L:\QQ) \otimes e^q)^p$$ converging to $HF(L,L)$ where
$(\;)^p$ means the total degree $p$. Moreover it collapes after a
finite number of steps.
\end{theorem}
\begin{proof}
We only prove the last statement. For the holomorphic disc $\beta$
with $\mu(\beta) \geq n+2$, the boundary $\delta_\beta$ is a zero map
because of Theorem \ref{dim}.
In monotone case, there exists only finitely many
homotopy classes $\beta$ with $\mu(\beta) < n+2$.
So, spectral
sequence collapes at a finite step, say $r_0$. In fact, $r_0$ may
be taken as the smallest number
 which satisfies $$(2r_0 -2) \geq (n+2).$$
\end{proof}
\section{Orientation}\label{sec:orientation}
We consider an orientation on the moduli space of
holomorphic discs with Lagrangian boundary condition. It is well-known
that moduli space of $J$-holomorphic discs is not always
orientable. (For example, consider the Lagrangian submanifold
$\RR P^2 \subset \CP^2$. The moduli space of constant discs
with boundary in $\RR P^2$ is non-orientable.)

In this section, we recall how to orient the moduli space of
$J$-holomorphic discs with a given spin structure.
In section \ref{sec:spinstructure}, we analyze how the
change of spin-structure of a Lagrangian submanifold will result
in the change of orientation described in this section. 
In section \ref{sec:ori}, we introduce necessary orientation conventions and
formulae, which will be used for the explicit computation for the
case of the Clifford torus.
In section \ref{sec:torusori} we show that there exists a standard spin-structure
for the Clifford torus and we describe how it determines the orientation of
moduli spaces of holomorphic discs.

We first recall the following theorem about orientability of the
moduli space in \cite{fooo:lift00}
\begin{theorem}[\cite{fooo:lift00}Theorem 21.1]
The moduli space of $J$-holomorphic discs is orientable, if $L
\subset(M,\omega)$ is a (relatively) spin Lagrangian submanifold.
Furthermore the choice of (relative) spin structure on $L$
determines an orientation on $\CM(L,\beta)$ canonically for all
$\beta \in \pi_2(P,L)$.
\end{theorem}
\begin{remark}
For simplicity, we will sketch the proof only when $L$ is a spin
manifold. For the relative spin case, see \cite{fooo:lift00}.
\end{remark}
\begin{proof} It suffices to show that the index of the linearized
operator is oriented. The linearized operator
$D\overline{\partial}$ for the $J$-holomorphic curve equation is
a first order elliptic differential operator with the same
symbol as the Dolbeault operator: for the $J$-holomorphic map
$w:(D^2,\partial D^2) \to (M,L)$ with $(p>2)$,
$$ D\overline{\partial}_w : W^{1,p}(D^2,\partial D^2;w^*TM, (w|_{\partial D^2})^*TL)
\to L^p(D^2;w^*TM). $$
It suffices to show that the index of the linearized
operator is oriented.
Since the zero order term does not affect
the index problem, we assume that the operator is the Dolbeault
operator $\overline{\partial}_w$
$$\overline{\partial}_{(w^*TM,(w|_{\partial D^2})^*TL)}:
W^{1,p}(D^2,\partial D^2;w^*TM, (w|_{\partial D^2})^*TL)
\to L^p(D^2;w^*TM).$$
 We recall how to determine a pointwise orientation of the index bundle from
\cite{fooo:lift00}
\begin{prop}[\cite{fooo:lift00} Proposition 21.3]\label{orientspin}
Let $E$ be a complex vector bundle over a disc $D^2$. Let $F$ be a
totally real subbundle of $E|_{\partial D^2}$ over $\partial D^2$.
We denote by $\overline{\partial}_{(E,F)}$ the Dolbeault operator
on $D^2$ with coefficient $(E,F)$,
 $$\overline{\partial}_{(E,F)} : W^{1,p} (D^2, \partial D^2 ; E, F)
 \to L^p(D^2; E)$$
Assume F is trivial and take a trivialization of F over $\partial D^2$. Then
the trivialization gives an orientation of the virtual vector space
 $Ker \, \overline{\partial}_{(E,F)} - Coker \,
\overline{\partial}_{(E,F)}$
\end{prop}
\begin{proof}
Here is a proof of the Proposition given in \cite{fooo:lift00}.
Suppose that the operators are surjective. (Otherwise we consider a
quotient of the target by a finite dimensional complex
subspace). 
By deforming the Hermitian
connection, we may assume that the totally real subbundle F is
trivially flat and the connection is product in a collar
neighborhood of $\partial D^2$. Let $C$ be a concentric circle in
the collar neighborhood of $\partial D^2$. If we pinch $C$ to a
point, we have the union of a disc $D^2$ and a 2-sphere $\CP^1$
with the center $O\in D^2$ and $S \in \CP^1$ identified. By the
parallel translation along radials, the trivial vector bundle F
extends up to C and its complexification gives a trivialization of
$E|_C$. Thus the bundle descends to $D^2 \cup \CP^1$. We
also denote this vector bundle by $E$. Then one can show that the
indices of the following two operators are isomorphic to each
other with the argument given in \cite{mcduff:jcqc94} Appendix A.
$$\overline{\partial}_{(E,F)} : W^{1,p} (D^2, \partial D^2 ; E, F)
 \to L^p(D^2; E \otimes T_{0,1}^* D^2)$$
$$\overline{\partial} : \{ (\xi_0, \xi_1) \in
W^{1,p} (D^2, \partial D^2 ; E, F) \times W^{1,p}(\CP^1,E) \,\, \vert
\, \,\xi_0(O) = \xi_1(S)\}$$
$$ \to L^p(D^2; E \otimes T_{0,1}^* D^2) \times
L^p(\CP^1; E \otimes T_{0,1}^* \CP^1)$$

If we have an element of the second index bundle, then use a cut-off
function to define an approximate element in the kernel of the
first operator. Then it projects onto the kernel of the first
operator. Hence we will get an orientation preserving isomorphism.

Since the real vector bundle F is trivialized, and by the above construction,
the kernel of the second operator is the kernel of the homomorphism:
\begin{equation}\label{split}
(\xi_0, \xi_1) \in  Hol (D^2,\partial D^2 : \CC^n,\RR^n)
\times Hol(\CP^1,E) \to  \xi_0(O) - \xi_1(S) \in \CC^n \cong E_S
\end{equation}
 Note that the kernel can be oriented by the orientation of
$\RR^n \cong Hol (D^2,\partial D^2 : \CC^n,\RR^n)$ since $Hol(\CP^1,E)$, and  $\CC^n$ carries a complex
orientation. This proves the Proposition. 
\end{proof}


Hence, we set $E =w^*TM, F = (w|_{\partial D^2})^*TL$ and apply 
Proposition \ref{orientspin} to determine the pointwise orientation
of the index bundle of $\overline{\partial}_{(w^*TM,(w|_{\partial D^2})^*TL)}$.

After fixing an orientation at one disc, say $w_0$, we can extend
this orientation to any disc $w$ in the path component of the
moduli space containing $w_0$. We consider a path
$w_t:(D^2,\partial D^2) \to (M, L)$ for $t \in [0,1]$ starting
from $w_0$ and ending at $w$. Since [0,1] is contractible, we have
a trivialization of $(w_t|_{\partial D^2})^*TL$, which gives an
orientation for $Index (\overline{\partial} w_t)$ for all $t\in
[0,1]$, i.e the orientation for $w$.

 Now under the assumption that $L$ is spin, we can show
that this assignment of orientation described above is independent
of a choice of paths. If there is a loop of holomorphic discs
$$w_\theta : (D^2 \times S^1, \partial D^2 \times S^1) \to (M, L),$$
$(w_\theta|_{\partial D^2}) ^* TL$ will be (stably) trivial over
$\partial D^2 \times S^1$ because $L$ is spin. So, it implies that
we will get a consistent orientation. This finishes the proof.
\end{proof}

\section{The changes of spin structures}\label{sec:spinstructure}
In this section, we analyze how the change of spin structures
affects the orientation of the moduli space.
 First we recall the
definition of a spin-structure given by Milnor \cite{milnor:ssm63}
\begin{definition}\label{spin}
A \textbf{spin structure} of an oriented vector bundle $E$ over
$X$ is a homotopy class of a trivialization of $E$ 
 over the 1-skeleton of $X$ which can be extended to
the 2-skeleton of $X$.
\end{definition}
\begin{remark} Note that orientation is a homotopy class of a
trivialization over the 0-skeleton which can be extended to
the 1-skeleton.
\end{remark}
The above definition is equivalent to the usual definition of 
the spin structure (for example, the definition in \cite{lawson:sg86}).
Recall that an oriented vector bundle $E$ over a manifold is called {\bf spin} if
its second Stiefel-Whitney class of $E$ is zero.
Here are basic properties of spin-ness
\begin{theorem}[\cite{lawson:sg86} Theorem 2.1.3]
Let $E$ be an oriented vector bundle of rank $\geq 3$ over $X$.
Then $E$ is spin if and only if for any compact surface $\Sigma$
and any continuous map $f:\Sigma \to X$, the bundle $f^*E$ is
trivial.
 Furthermore, if $E$ is spin,
then the distinct spin structures on E are in one to one
correspondence with the elements of $H^1(M;\ZZ/2\ZZ)$.
\end{theorem}

 In case the rank of the bundle $E$ is two or less, we add trivial
vector bundle to the bundle $E$ , and we will get a stable
trivialization instead of a trivialization by the previous Theorem.
But stable trivialization is good enough to deal with the orientation 
problem

Let $w:(D^2,\partial D^2) \to (M,L)$ be a holomorphic disc with
Lagrangian boundary condition. Recall that orientation of the
moduli space of holomorphic discs is determined by the
trivialization of $(w|_{\partial D^2})^*TL$. We will prove the
following Theorem which is crucial to understand the change of
orientation.
\begin{theorem}\label{sgnchange}
If we reverse the homotopy class of a trivialization of
$(w|_{\partial D^2})^* TL$, then orientation given on the index
bundle of $\overline{\partial}_{w^*TM,(w|_{\partial D^2})*TL}$ will be
reversed.
\end{theorem}
\begin{remark}
Note that for any orientable vector bundle over $S^1$, there exist only two 
homotopy classes of a (stable)
trivializations.
The above theorem was already stated without proof as a remark 21.6 of
\cite{fooo:lift00} with an example.
\end{remark}

With this Theorem \ref{sgnchange} in hand, one can analyze the change of
orientation as follows:
First we fix a spin structure of a Lagrangian submanifold $L$.
Hence, for any holomorphic disc $w:(D^2,\partial D^2) \to
(M,L)$, this determines the homotopy type of the trivialization of
$(w|_{\partial D^2})^* TL$. By \cite{fooo:lift00} Proposition 21.3, it
determines the orientation of the moduli space of holomorphic
discs $\widetilde{\CM}(\beta)$.
To analyze the orientation of $\widetilde{\CM}(\beta)$ for a {\em
different spin structure}, we note that a change of a spin
structure will result in a change of homotopy type of a trivialization 
of $TL$ where
both change correspond to $H^1(L;\ZZ/2\ZZ)$. Hence for any
holomorphic disc $w$, if the homotopy class of a trivialization of
$(w|_{\partial D^2})^* TL$ is reversed due to the change of the
spin structure, then the induced orientation for
$\widetilde{\CM}(\beta)$ will be reversed. This will be exactly
the way we will calculate Floer cohomology of the Clifford torus
with various spin-structures.

To prove the Theorem \ref{sgnchange}, we need an index theorem for the
holomorphic discs proved by Silva \cite{silva:psfh97}.
We will state it briefly here.
\begin{definition}
A bundle pair $(T,\lambda)$ is a complex vector bundle $T$ over
$D^2$ and
 a real vector bundle $\lambda$ over $\partial D^2$ such that
$\lambda \otimes \CC$ is identified with $T|_{\partial D^2 }$.
\end{definition}
For such a pair, we define $ind(T,\lambda)$ as follows. To
incorporate boundary conditions, we restrict the domain of
$\overline{\partial}$ :
$$\Gamma_\lambda(T):=\{ s \in \Gamma(T) | s:\partial D^2 \to \lambda
 \subset T|_{\partial D^2} \}$$
then, a canonical Cauchy-Riemann operator with addtional boundary
condition gives an elliptic boundary value problem.
$$\overline{\partial} :\Gamma_\lambda(T) \to \Gamma(T \otimes T^{0,1}D^2)$$
We define $ind(T,\lambda)$ to be the index bundle of this
operater. It depends only on the homotopy type of $(T,\lambda)$
and is additive
under taking direct sums.

Now we consider a family of discs, parametrized by a compact space
$X$.
\begin{definition}
A bundle data $(T,\lambda)$ is a unitary bundle $T$ over $D^2
\times X$ and
 an orthogonal bundle $\lambda$ over $\partial D^2 \times X$ such that
$\lambda \otimes \CC$ is identified with $T|_{\partial D^2 \times
X}$.
\end{definition}
By choosing a continuous family $\{\overline{\partial}_x\}_{x\in
X}$ of Cauchy-Riemann operators in the fibres, we obtain an index
bundle
$$ind(T,\lambda) \in KO(X)$$
Again, the index depends only on the homotopy type of
$(T,\lambda)$ and is additive under taking direct sums.

For a bundle data $(T,\lambda)$, assume $T$ is trivial of rank
$n$. Then by fixing a trivialization of $T$ over $D^2 \times X$, we can
specify $\lambda$ as a map $\phi_\lambda : \partial D^2 \times X
\to U(n)/O(n)$. If $n$ is large compared to the dimension of $X$, we
can replace $U(n)/O(n)$ by its stable limit $U/O$. Let $x_0$ be a
basepoint for $X$ and also assume that $\phi_\lambda$ is constant
on $\partial D^2 \times \{x_0\}\cup\{1\}\times X$. The index
bundle then necessarily has rank $n$ since it is trivial over $x_0$.
Subracting $n$ from it, we obtain an element of $\widetilde{KO}(X)
\cong [X,BO]$.
\begin{theorem}[Silva~\cite{silva:psfh97}]\label{silvaindex}
The above construction gives an isomorphism of abelian
groups $$ind:[\Sigma X, U/O]_* \to [X,BO]$$ where the addition is
given by taking direct sum. 
Here, $\Sigma X$ denotes a reduced suspension of $X$. $[,]_*$ denotes
homotopy class of based maps $\;\;\;\Box$
\end{theorem}
\begin{remark}
Note that the left hand side is only defined under the assumption that
the bundle $T$ is trivial and $\phi_\lambda$
is constant on  $\partial D^2 \times \{x_0\}\cup\{1\}\times X$. 
Hence one can not define this isomorphism  for every bundle pair.
But one can take the direct sum of the given vector bundle with 
a certain bundle pair to 
define such a map ( for details, see \cite{silva:psfh97}).
\end{remark}

We will apply this index Theorem for the case $X=S^1$.
Then we have 
$$[\Sigma S^1,U/O]_* = [S^2, U/O]_* = \pi_2(U/O)
\underbrace{\cong}_{Ind} \widetilde{KO}(S^1) \cong \pi_0(O)\cong
\ZZ/2\ZZ$$ So a non-trivial generator of $\pi_2(U/O)$ will give rise
to
the non-orientable index bundle over $S^1$.
But in view of the homotopy exact sequence
$$\pi_2(U) \to \pi_2(U/O) \to \pi_1(O)$$
with $\pi_2(U) \cong 0$, the non-trivial element
of $\pi_2(U/O)$ corresponds to the non-trivial element of $\pi_1(O)$.
This correspondence is the main reason for the Theorem
\ref{sgnchange}: The
change of the homotopy class of a trivialization of the real
vector bundle $\lambda$, which in fact comes from the twisting
of a trivialization by non-trivial element of $\pi_1(O)$, will
reverse the orientation of the index bundle
given by Proposition 26.2 \cite{fooo:lift00}.

\begin{proof}
We will construct a bundle pair which contains both homotopy classes of
trivializations and we will show that its index bundle is
non-orientable using Theorem \ref{silvaindex}.
We start with the case of Maslov index 0.

Consider the trivial bundle $(D^2 \times [0,1]) \times \CC^N$ over 
$(D^2 \times [0,1])$.
On the base, by identifying $D^2\times\{0\}$, $D^2\times\{1\}$, we get $D^2
\times S^1$.
And we glue the fibers $\partial D^2\times\{0\}\times \RR^n$,
$\partial D^2\times\{1\}\times \RR^n$ by homotopically non-trivial
loop
$\gamma :\partial D^2 \to SO(n)$ with $\gamma(1) = Id \in SO(n)$.
$i.e.$ for $z \in D^2, x \in \RR^n$, we identify
$$(z,0,x) \sim (z,1,\G(z)x)$$
The inclusion of
    $\pi_1(SO(n)) \to \pi_1(SU(n))$ is trivial
since
$\pi_1(SU(n))\cong 0$. So we can extend the map $\gamma$ to
$\Gamma : D^2 \to SU(n)$. Then we identify $D^2\times\{1\}\times
\CC^n$ with $D^2\times\{0\}\times \CC^n$ by the map $\Gamma$. Note
that this identification matches with the one given on
$\RR^n$ before. We denote the resulting bundle data as
$(T,\lambda)$.
We can give a trivialization of the bundle $T$ as follows:
Let $C_{Id}: D^2 \to U(n)$ be the constant map $C_{Id}(z) = Id \in
U(n)$
for $z\in D^2$.

First, there is a homotopy $H:D^2 \times [0,1] \to U(n)$ between the
two maps $C_I$ and $\Gamma$ such that for $ z\in D^2, t \in [0,1]$
\begin{equation}
\begin{cases}
H(z,0) = C_{Id}(z)\\
H(z,1) = \Gamma(z)\\
H(1,t) \equiv Id \in U(n)
\end{cases}
\end{equation}
Then, we define
$$\Psi : D^2\times [0,1] \times \CC^n \to D^2\times [0,1] \times \CC^n$$
$$ (z,t,x) \to (z,t,H(t,z)x)$$
This map $\Psi$ defines a trivialization of $T$ : We identified
$(z,1,x)$ with $(z,0,\Gamma(z)(x))$ and
$$\Psi(z,1,x) = (z,1,H(1,z)x) = (z,1,\Gamma(z)x)$$
$$\Psi(z,0,\Gamma(z)(x)) =(z,0,H(0,z)\Gamma(z)(x)) = (z,1,\Gamma(z)x)$$
Hence, the trivialization $\Psi:T \to D^2 \times S^1 \times \CC^n$
is well-defined where $S^1$ is given by $\RR/\ZZ$.
Under the trivialization $\Psi$, we define a map 
$\phi:\partial D^2 \times S^1 \to U(n)/O(n)$ as follows: For 
$z \in \partial D^2, t \in S^1$
$$ (z,t) \mapsto [\lambda_{(z,t)}] \in U(n)/O(n)$$
where $[\lambda_{(z,t)}]$ is an element in the Lagrangian Grassmaniann
corresponding to the Lagrangian subspace $\lambda_{(z,t)} \subset \Psi(T_{(z,t)})
\cong \CC^n$.

Then from the construction of $\Psi$, we have
\begin{equation}
\phi|_{\partial D^2 \times 0} \equiv 
\phi|_{\partial D^2 \times 1} 
\equiv Id \in  U(n)/O(n) 
\end{equation}
$$\phi|_{\{1\}\times S^1} \equiv Id$$
By moding out $\partial D^2 \times \{0\} \cup \{1\}\times S^1$ from
$\partial D^2 \times S^1$, we may consider $\phi$ to be a map
from $S^2$ to $U(n)/O(n)$.

It is not hard to see that $\phi$ gives an non-trivial element in $\pi_2(U(n)/O(n))$.
So the bundle data $(T,\lambda)$ has non-orientable Cauchy Riemann
index bundle by the Silva's index Theorem. And it implies that
orientations given by these two different homotopy classes of trivializations 
can not be same.

Now we study the case of Maslov index 2.
Basically we will use the same technique as Maslov index 0 case.
In $\CC^n$, for $z \in \partial D^2$, our Lagrangian subspaces
will be $z\cdot \RR \times \RR^{n-1}$. Take the trivial bundle
$(D^2 \times [0,1]) \times \CC^N$. By identifying
$D^2\times\{0\}$, $D^2\times\{1\}$, we get $D^2 \times S^1$.

We glue the Lagrangian subspaces of $\CC^n$ fibers, $\partial D^2 \times \{0\}\times z\cdot \RR
\times \RR^{n-1}$, $\partial D^2 \times \{1\}\times z\cdot \RR
\times \RR^{n-1} $ with non-trivial loop in $SO(n)$  as follows:
We define a map $R:S^1 \to U(n)$ as 
$R(e^{i\theta}) = \mbox{diag}(e^{-i\theta},1,\cdots,1) \in
U(n)$. Along $\partial D^2$, we identify two Lagrangian fibres by 
$R^{-1}\circ \gamma
\circ R$ where $\gamma$ is the non-trivial loop in $SO(n)$ used
in Maslov index 0 case. This identification can be extended to the
whole fiber $\CC^n$ and also over $D^2$, since $R^{-1}\circ \gamma \circ R$ 
is a loop in $SU(n)$.
Let $(\widetilde{T},\widetilde{\lambda})$ be the bundle data obtained by gluing each end
with this identification.

To show that the bundle data $(\widetilde{T},\widetilde{\lambda})$ has non-orientable
Cauchy-Riemann index bundle, we will show that $(\widetilde{T},\widetilde{\lambda})$
still gives the non-trivial element of
$\pi_2(U/O)$ after some modification.

First, we will make a direct sum $(\widetilde{T},\widetilde{\lambda})$ with a bundle pair
with Maslov index $(-2)$:
$$ i.e.\,\,(D^2\times I \times \CC, \partial D^2(z)
\times I \times z^{-1} \cdot \RR)$$ 
Denote the resulting 
bundle data as $(\widetilde{T}',\widetilde{\lambda}')$.
And when we glue two ends 0 and 1,
we glue this extra $\CC$-fiber without twisting. Since this extra
bundle pair has an orientable index bundle, and index bundle is
additive upon direct sums. So it is enough to show that
the bundle data $(\widetilde{T}',\widetilde{\lambda}')$
obtained this way is non-orientable.

Now we will find a trivialization of $\widetilde{T}'$ as Maslov index 0 case.
We will regard $\gamma(z) \in SU(n)$ and $\Gamma(z)\in SU(n)$ 
as in $SU(n+1)$ extending by 0 except the last diagonal entry
where we extend it by 1.

Now, let $$ F(e^{i\theta}) =  diag( e^{-i\theta},1,\cdots,1,e^{+i\theta})
 \in SU(n+1)$$

The map $F:\partial D^2 \to SU(n+1)$
can be extended over $D^2$ and we will denote the extention again by 
$F:D^2 \to SU(n+1)$.

Again let $C_{Id}: D^2 \to SU(n+1)$ be a constant map to $Id \in
SU(n+1)$. We get a homotopy $\widetilde{H}(z,t):D^2 \times [0,1] \to U(n+1)$
 between $C_{Id}$ and $F^{-1} \circ \gamma \circ F$ with
the following properties: For $z \in D^2, t\in S^1$,
\begin{equation}
\begin{cases}
\widetilde{H}(z,0) = C_{Id}(z) \equiv Id \in SU(n+1)\\
\widetilde{H}(z,1) = F^{-1} \circ \Gamma \circ F(z)\\
\widetilde{H}(1,t) = Id \in U(n+1)
\end{cases}
\end{equation}

Then, we define
$$\widetilde{\Psi} : D^2\times [0,1] \times \CC^{n+1} \to D^2\times [0,1] \times \CC^{n+1}$$
$$ (z,t,x) \to (z,t,F(z)\widetilde{H}(t,z)x)$$
One can check that this map $\widetilde{\Psi}$ gives a well-defined trivialization 
of the bundle $\widetilde{T}'$ as before.
Note that since we composed $F(z)$ in the trivialization $\widetilde{\Psi}$, 
the real vector bundle $\widetilde{\lambda}'$
will map
to constant Lagrangian subspace $\RR^n \subset \CC^n$ over $\partial D^2 \times 0$.

Under this trivialization, we can similarly define a map by considering
the image the real bundle $\widetilde{\lambda}'$ in the Lagrangian Grassmaniann,
$$\phi :\partial D^2 \times S^1 \to U(n+1)/O(n+1)$$
Then, 
$$\phi|_{\partial D^2 \times 0} 
 \equiv 
\phi|_{\partial D^2 \times 1} 
\equiv Id \in  U(n+1)/O(n+1)$$ 
$$\phi|_{\{1\}\times S^1} \equiv Id$$
By moding out $\partial D^2 \times \{0\} \cup \{1\}\times S^1$ from
$\partial D^2 \times S^1$, we may consider $\phi$ to be map
from $S^2$ to $U(n)/O(n)$.

As before, it defines a nontrivial element in $\pi_2(U/O)$.
Hence the bundle pair $(\widetilde{T}',\widetilde{\lambda}') $ has non-orientable Cauchy
Riemann index bundle by Silva's index Theorem, which implies
that original bundle pair $(\widetilde{T},\widetilde{\lambda}) $ also has
non-orientable index bundle.

Other cases can be done similarly. 
This finishes the proof of the Theorem \ref{sgnchange}
\end{proof}

\section{Orientation conventions and formulae}\label{sec:ori} 

In this section, we will fix some basic
conventions concerning orientations. These conventions agree with
the ones defined in \cite{fooo:lift00}.

\begin{itemize}
\item We will assume that all circles are oriented counter-clockwise.

\item We orient the Clifford torus $T^n$ as
a torus $(S^1)^n \subset U_0 \cong \CC^n$ where $U_0 = \{ z_0 \neq 0 \} \subset \CP^n$.
where each $S^1 \subset \RR^2 \cong \CC$ is oriented counter-clockwise.

\item For a Clifford torus, we have a torus action $(S^1)^n$ action on it
given by 
$$(e^{i\theta_1},\cdots,e^{i\theta_n})\cdot [z_0;\cdots;z_n]
\mapsto [z_0;e^{i\theta_1}z_1;\cdots;e^{i\theta_n}z_n].$$

\item The elements of $PSL(2:\RR)$ can be written as
$e^{i\theta}(\frac{z-\alpha}{1-\overline{\alpha}z})$ for
$\alpha \in D^2$. We will orient $PSL(2:\RR)$ as $S^1 \times D^2$ where
the latter carries a complex orientation.

\item Let $X$ be an oriented smooth manifold with boundary $\partial X$.
Then we define an orientation on $\partial X$ so that 
$$ T_*X \cong \RR_{out} \times T_*(\partial X).$$
is an isomorphism of oriented vector spaces. Here $\RR_{out}$ is an $\RR$
oriented by outer normal vector.

\item Let $G$ be a Lie group given an orientation. When $G$ acts on an oriented
manifold $X$ smoothly and freely, then we define an orientation of
the quotient space $X/G$ so that 
$$T_*X \cong T_*(X/G) \times \textrm{Lie}\;G$$
is an isomorphism of oriented vector spaces. Here Lie $G$ is the Lie algebra
of $G$.

\item We orient the moduli space $\widetilde{\CM}_m(\beta)$ as 
$\widetilde{\CM}(\beta) \times (\partial D^2)^m$.

\item In \cite{fooo:lift00}, the orientation of the fibre product
$X_1 \,\times_Y X_2$ is given
for the case when the maps $f_i : X_i \to Y$ are submersions. 
Here we specify the orientation for the case that the map $f_2:X_2 \to Y$ 
is an embedding.
This will be used throughout our computation.

Let $X$, $L$ and $P$ be oriented smooth submanifolds and
let $f : X \to L$ be a submersion and $i : P \to L$ be an embedding.
Here we will regard $P$ as a submanifold of $L$.
By $x, l, p$ we denote the dimension of $X, L, P$.
Take a point $q \in f(X) \cap P$.
We can choose an oriented basis $< u_1, \dots, u_l> \in T_q L$  so that
$< u_1, \dots, u_p> \in T_q P$ becomes an oriented basis for the given
orientations of $L$ and $P$.
Since $f$ is a submersion, we can choose $< v_1, \dots, v_l> \in T_p X$ for
some $p \in f^{-1} (q)$ such that $(df)_p (v_k) = u_k$ for $k = 1,\dots,l$.
Then, we can choose a basis $< \eta_1, \dots, \eta_{x-l}> \in Ker (df_p)$
such that $< \eta_1, \dots, \eta_{x-l}, v_1, \dots, v_l, > $ is the
given orientation of $T_p X$. 
Then we define an orientation on the fibre product $X \,_f \times _i P$ so
that $<\eta_1, \dots, \eta_{x-l},u_1,\dots,u_p>$ becomes an oriented basis.
\end{itemize}
In this setup, it is easy to see that
\begin{lemma}
$$ \partial (X \,_f \times P) = \partial X \,_f\times P \bigsqcup
(-1)^{x-l} X \,_f \times \partial P $$
\begin{equation}
= \partial X \,_f \times P \bigsqcup
(-1)^{x+l} X \,_f \times \partial P \label{eq:ori1}
\end{equation}
\end{lemma}
\begin{proof}
Recall that oriented basis of $X$ was written as
$< \eta_1, \dots, \eta_{x-l}, v_1, \dots, v_l, > $. We write this as 
$[X] = [X^0] \times [L]$ where $[X^0]$ represents the basis 
$< \eta_1, \dots, \eta_{x-l}>$. 
From our convention for the boundary orientation, we have
$$[X] = [\RR_X] \times [\partial X].$$
Hence, we write $$[P] = [\RR_P] \times [\partial P].$$
The orientation of $X \times_L P$ can be written as 
$$[X \times_L P] = [X^0] \times [P].$$
Hence
\begin{eqnarray*}
[\RR_{X \times_L P}] \times [\partial(X \times_L P)] &=& [X \times_L P] =
[X^0] \times [P]\\
&=& [\RR_X]\times [(\partial X)^0] \times [P] \bigsqcup [X^0] 
\times [\RR_P] \times [\partial P] \\
&=& [\RR_X]\times [(\partial X)^0] \times [P] \bigsqcup 
(-1)^{dim X^0}[\RR_P]\times [X^0] \times [\partial P] \\
&=& [\RR_ X] \times [(\partial X) \times_L P]
\bigsqcup (-1)^{x-l} [\RR_P]\times [X \times_L \partial P].
\end{eqnarray*}
\end{proof}

Another formula we use in the proof of Theorem \ref{squarezero} is 
the orientation formula for the gluing from Proposition 23.2 in \cite{fooo:lift00}.
\begin{prop}[\cite{fooo:lift00} Proposition 23.2]
\begin{equation}
\partial \CM_2(A+B) = (-1)^{dim L + 1} \CM_2 (A)_{ev_1} \times_{ev_0}
\CM_2(B). \label{eq:ori2}
\end{equation}
\end{prop}

\section{Orientation of the moduli space for the Clifford torus}\label{sec:torusori}

In this section, we show that there exists a natural spin structure of
the Clifford torus, which we denote by {\em standard spin structure}.
Under the standard spin structure, it is rather easy
to determine the orientation of the moduli space as 
described in Theorem \ref{orientspin}.
\begin{prop}\label{prop:standard}
There exists a {\em standard} spin structure of the Clifford torus. 
Or equivalently,
there exists a natural homotopy class of a trivialization of the tangent bundle $T(T^n)$
of the Clifford torus.
\end{prop}

\begin{proof} Let $S^1:= e^{i\theta}$  be the unit circle embedded in $\CC$. The
tangent bundle of $S^1$ has a natural trivialization given by $S^1
\times \RR \cdot \PD{\theta}$. Similarly there is a natural
trivialization of the tangent bundle of $(S^1)^n \subset \CC^n$.
The Clifford torus $T^n$ sits inside the intersection of 
$n+1$ standard open covers
$U_i(\cong \CC^n)=\{z_i \neq 0\} \subset \CP^n$. 
So, each open cover induces a trivialization
of tangent bundle of $T^n$. One can check that the 
trivializations of $T(T^n)$ obtained with each open set $U_i$ 
are in the same homotopy class:
Because the transition matrices between these trivializations
are constant matrices, which implies that there is
no twisting of frames. 
By permuting coordinates to have positive determinants, if necessary, 
the trivializations
 induced in each open set are 
in the same homotopy class.
This is what we mean by the standard spin structure of $T^n$.
\end{proof}

But we need to fix a trivialization in this homotopy class to fix an orientation.
We will fix the
trivialization to be the one obtained
from the open set $U_0 \subset \CP^n$.

Now, we discuss the orientation of the moduli space of 
holomorphic discs with boundary on $T^n$. The discussion is based on
the classification Theorem of such holomorphic discs in section
\ref{sec:classify}.
Let $\widetilde{\CM}(\beta)$ be the space of holomorphic discs
representing the homotopy class $\beta \in \pi_2(M,L)$ as defined in 
Definition \ref{modulinoquo}.
The orientation of $\widetilde{\CM}(\beta)$ can be determined by
the Proposition \ref{orientspin} after we fix the spin structure.
We start with an example.

For a homotopy class $\beta_0 \in \pi_2(\CP^n,T^n)$, we will see that
the moduli space $\widetilde{\CM}(\beta)$ is 
\begin{equation}\label{example}
\{ \;[ \frac{z-\alpha}{1-\overline{\alpha}z}:e^{i\theta_1}:\cdots:
e^{i\theta_n}] |
\alpha \in int\, D^2, \theta_i \in S^1 \}
\end{equation}
Since $\alpha \in D^2 \subset \CC$ carries a complex orientation, 
the orientation of $\widetilde{\CM}(\beta)$ is determined by the orientation of 
$(\theta_1,\cdots,\theta_n) \in (S^1)^n$.
With the standard spin structure, it will oriented as 
$(\PD{\theta_1},\cdots,\PD{\theta_n})$.

From now on we fix the standard spin structure.
Let $w:(D^2,\partial D^2) \to (\CP^n,T^n)$ be a holomorphic disc.
Recall that in Proposition \ref{orientspin}
we had a decomposition of the tangent space of
$\widetilde{\CM}(\beta)$  as a kernel of the homomorphism
$$(\xi_0, \xi_1) \in  Hol (D^2,\partial D : \CC^n,\RR^n)
\times Hol(\CP^1,E) \to  \xi_0(O) - \xi_1(S) \in \CC^n.$$
Here $Hol (D^2,\partial D : \CC^n,\RR^n)$ is in fact just $\RR^n$
and this $\RR^n$ comes from the trivialization $T(T^n)$ along $w|_{\partial D^2}$.
It is not hard to see that this $Hol (D^2,\partial D : \CC^n,\RR^n)$
corresponds to the subspace of tangent space $T_w(\widetilde{\CM}(\beta))$
which is given by the translation of disc $w$ along the tangent directions
of the Lagrangian submanifold $T^n$ under the standard spin structure.
Therefore, $Hol (D^2,\partial D : \CC^n,\RR^n)$ is oriented by our choice of 
natural trivialization in the previous Proposition. 
All other factors in the above decomposition carries complex orientations. 
This gives the orientation of $\widetilde{\CM}(\beta)$. 

Under non-standard spin structures, such a direct analysis is not possible, but
we can still assign the orientation of the moduli space as described in 
the paragraph after Theorem \ref{sgnchange}.\\

Now we compute the orientation of $\CM_1(\beta_i)$ with our orientation convention.
Recall that 
we orient the moduli space with marked points $\widetilde{\CM}_m(\beta)$ as
$[\widetilde{\CM}(\beta)] \times [\partial D^2]^m$.
where $[,]$ means the oriented basis of the tangent space.
The moduli space
$\CM_m^{reg}(\beta)$ is oriented as $[\widetilde{\CM}(\beta)] \times
[(\partial D^2)^m] / PSL(2:\RR).$ or
$$[\widetilde{\CM}(\beta)] \times [(\partial D^2)^m] =
[\CM_n^{reg}(\beta)]\times [PSL(2:\RR)]$$ where $[PSL(2:\RR)]$
represents a frame at the tangent space of each disc in $\widetilde{\CM}_m(\beta)$
which is
given by $[PSL(2:\RR)]$ action on $[\widetilde{\CM}(\beta)] \times
[(\partial D^2)^n]$.
If we only consider the holomorhic discs with Maslov index 2, 
Then the homotopy classes $\beta_i$'s are minimal, hence 
$\CM_1^{reg}(\beta_i) = \CM_1(\beta_i).$
Hence, we have $$\widetilde{\CM}(\beta_i)/PSL(2:\RR) \cong \CM_1(\beta_i).$$
\begin{prop}\label{tn}
Let $\beta_i \in \pi_2(\CP^n,T^n)$ be the homotopy class described in Proposition
\ref{torusmonotone} for $i =0,1,\cdots,n$.
Then the evaluation map $ev_0 : \CM_1(\beta_i) \to T^n$ is an orientation 
preserving homeomorphism for all $i=0,1,\cdots,n$.
\end{prop}
\begin{proof}
It is easy to see that $ev_0$ is a homeomorphism by the classification
Theorem \ref{classify}. 
To find out the orientation of $\CM_1(\beta)$, we need to specify the orientation
of $\widetilde{\CM}(\beta_i)$ and $[PSL(2:\RR)]$.
Since $[\widetilde{\CM}_1(\beta)]\cong [\widetilde{\CM}(\beta)] \times [\partial
D_0^2]$, we have 
$$[\CM_1(\beta_i)] =([\widetilde{\CM}(\beta)] \times
[(\partial D_0^2)]) / PSL(2:\RR).$$
Recall that $T(\widetilde{\CM}(\beta_i))$ have a decomposition as 
$[T^n] \times [D^2]$ (see the expression \ref{example}).
Hence, by taking a quotient of $D^2 \subset PSL(2:\RR)$ which carries a complex orientation,
$$[\CM_1(\beta_i)] = ([T^n] \times [D^2] \times [\partial D_0^2])/[S^1] \times 
[D^2]$$
$$= ([T^n] \times [\partial D_0^2])/[S^1] $$
Here  an element $e^{i\theta} \in  S^1 \subset PSL(2:\RR)$ acts on a holomorphic disc 
$w$ as $e^{i\theta} \cdot w(z) = w (e^{-i\theta}z)$
and it acts on a marked point $z_0$ as $e^{i\theta}z_0$.
So under the evaluation map $ev_0$ ,
we obtain $$[(\CM_1(\beta_i),ev_0)] \cong [T^n].$$
This finishes the proof.
\end{proof}
\section{Maslov index formula for discs in $\CP^n$ with boundary in $T^n$}

The most challenging part in computing the Floer cohomologies is to 
classify all the holomorphic($J$-holomorhpic) discs with Maslov index $\leq$
$n+1$.
In the case of the Clifford torus, or generally for torus fiber in compact toric
manifolds, the following formula is a fundermental tool.
Later on, by using this index formula, 
we will classify all holomorphic discs 
with boundary on the Clifford torus for any Maslov indices.

\begin{theorem}\label{indexformula}
For a holomorphic disc  $w:(D^2,\partial D^2) \to (\CP^n,T^n)$,
the Maslov index of the disc $w$ is twice the sum of intersection
multiplicities between the image of the disc $w$ with hyperplanes $H_i$ for
$i= 1 \dots n$, where $H_i$'s are hyperplanes defined by $z_i=0$
in $\CP^n$
\end{theorem}
\begin{proof}
We first prove the following elementary lemma regarding the Maslov index 
of a map.
\begin{lemma}\label{cot}
Let $L$ be a Lagrangian submanifold whose tangent bundle $TL$ is
trivial. Let $o_L$ be the zero section of the cotangent bundle of
$T^*L$.
Let $\Sigma$ be a smooth Riemann surface with boundary.
Then, for any smooth map $w:(\Sigma,\partial \Sigma) \to
(T^*L, o_L)$, the Maslov index of $w$ is zero.
\end{lemma}
\begin{proof} First, consider the case that the image of $w(\Sigma)$
is entirely contained in $o_L$. At the zero section $o_L$ of the
cotangent bundle, there exists a canonical splitting of $T(T^*L)
\cong TL \oplus T^*L$. From the trivialization of $TL$, we obtain
a trivialization of $TM|_{o_L} \cong TL \oplus J_0 (TL)$. Hence
for the pull back the above trivialization by $w$, the Maslov
index of the map $w$ is zero.

For general cases, we can homotope $w$ to $pr(w)$ where $pr:T^*L
\to o_L$ is the projection map of the cotangent bundle of $L$, and
we may replace the homotopy by a smooth one. Hence then $\mu(w) =
\mu(pr(w)) = 0$.
\end{proof}

Note that $\CP^n \setminus (H_0 \cup H_1 \cup \dots \cup H_n)$
can be identified with the cotangent bundle of $T^n$, which will be used
crucially later in the proof.

First, consider the case that a disc $w$ does not meet any 
hyperplanes $H_i$s at all. 
 Then the disc is in fact in the cotangent bundle of $T^n$. 
From the Lemma \ref{cot}, its Maslov index is zero, hence the 
Theorem holds in this case.

To consider the general discs, we first write the map in terms of 
the homogeneous coordinate functions.
\begin{lemma}
For a holomorphic disc $w:(D^2,\partial D^2) \to (\CP^n, T^n)$, we
can write the map as $$[\G_0(z):\G_1(z): \cdots :\G_n(z)]$$
where $\G_i(z):D^2 \to \CC$ is a holomorphic function for 
$i=0,\cdots,n$ with $\cap_{j=0}^n Zero(\G_j) = \phi $.
\end{lemma}
\begin{proof} There is the holomorphic line bundle $\mathcal{O}(1)$ over $\CP^n$
whose global sections are generated by $z_0,z_1,\cdots,z_n$.
Now consider the pull-back bundle $w^*\mathcal{O}(1)$ and 
over the disc, we fix its holomorphic trivialization $\Psi :
w^*\mathcal{O}(1) \to D^2 \times \CC$. 
Let $\G_i(z) = \Psi(w^*z_i)$. 
\end{proof}

Now we assume that there exists at least one intersection between 
the image of the map $w$ and the given hyperplanes, where 
one of $\G_i(z)$ becomes zero.
We label by $p_1, p_2, \cdots p_m \in D^2$ every point where one of the $\G_i(z)$ becomes zero. 
We find disjoint open balls $B_i(\epsilon) \subset D^2$ centered at $p_i$
with fixed radius $\epsilon$ for sufficiently small $\epsilon$ 
for all $i=1,2,\cdots,m$.
Our stratergy is to deform the map $w$ inside the ball $B_i$ 
so that we can decompose the disc into 
several regions whose boundary satisfies Lagrangian boundary
condition and then we will compute the Maslov index using the decomposition.

At $p_1$, we may assume without loss of generality that
$$\G_0(p_1)=\G_1(p_1)=\cdots=\G_s(p_1)=0$$ for $0 \leq s < n$.
Denote by $d_i$ the order of zero of $\G_i$ at $p_1$ for 
$0 \leq i \leq n$. Then, $d_i >0$ for $0 \leq i \leq s$ and 
$d_i =0$ for $s < i \leq n$.
We may further assume that $p_1 =0 \in D^2$ for simplicity.
Recall that $B_1(\epsilon)\subset D^2$ is a ball centered at 0 of radius
$\epsilon$. 
Since, $\G_i$ does not have common zero, $\G_n(B_1)$ is away from zero.
Hence the image of the map $w$ is contained in the open set $U_n$.

Define $f_i:B_1 \to U_n(\cong \CC^n)$ for $i = 0,1,\cdots,n-1$ as
$$f_i(z) = \frac{\G_i(z)}{\G_n(z)}.$$
These $f_i$'s are holomorphic functions.
For $i = 0,1,\cdots,n-1$,
we can choose $a_i \in \CC$ and holomorphic functions $R_i(z):B_1 \to \CC^n$
with
$$f_i(z) = a_i z^{d_i} + R_i(z)$$ 
where $R_i(z) = O(|z|^{d_i+1})$.

Basically, we want to deform the map $f_i(z) = a_i z^{d_i} + R_i(z)$ to
$(a_i z^{d_i}/|a_i|(\frac{\epsilon}{2})^{d_i} + 0)$ inside the ball 
$B_1(\frac{\epsilon}{2})$ without changing the map near the boundary of 
the ball $B_1(\epsilon)$. 
Note that the constants are chosen to map $\partial B_1(\frac{\epsilon}{2})$
to the Clifford torus.

Here are two cut-off type functions $\xi_c,\eta: \RR \to \RR$.
$$ \xi_c (x) = \left\{ \begin{array}{ll}
   1 & \textrm{if}\,\,  |x| \geq \frac{2 \epsilon }{3} \\
   \frac{1}{c} & \textrm{if}\,\, |x| \leq \frac{\epsilon}{2}
   \end{array} \right.$$
$$ \eta(x) =  \left\{ \begin{array}{ll}
    1 & \textrm{if}\,\, |x| \geq \frac{2 \epsilon }{3} \\
   0  & \textrm{if} \,\, |x| \leq  \frac{\epsilon}{2}
\end{array} \right.$$
We extend $\xi_c,\eta$ smoothly over $\RR$ with 
$\frac{1}{c} \leq \xi_c \leq 1, 0 \leq \eta \leq 1$.

We also define the deformation between $\xi_c,\eta$ and the constant
function $1$ as
$$\xi_c^t(x) = (1-t)\cdot 1 + t \xi_c(x)$$
$$\eta^t(x) = (1-t) \cdot 1 + t \eta(x)$$

Now we choose the constants $c_i = |a_i|(\frac{\epsilon}{2})^{d_i}$.
Let 
$$f_i^t(z) = \xi_{c_i}^t(|z|)a_i z^{d_i} + \eta^t(|z|) R_i(z)$$
Then, $f_i^0(z) = f_i(z)$ and $f_i^t$ gives the smooth deformation of the 
original map $w$ to a new map, say $w_1:(D^2,\partial D^2) \to (\CP^n,T^n)$, 
where $w_1|_{B_1(\frac{\epsilon}{2})}$ can be written as

\begin{equation}\label{map}
[ \frac{a_0 z^{d_0}}{|a_0|(\frac{\epsilon}{2})^{d_0}}:\cdots:
\frac{a_s z^{d_s}}{|a_s|(\frac{\epsilon}{2})^{d_s}}:\frac{a_{s+1}}{|a_{s+1}|}:
\cdots:\frac{a_{n-1}}{|a_{n-1}|}:1].
\end{equation}

We perform the same deformation for $p_2,p_3,\cdots,p_m$ inside the ball
$B_2,\cdots,B_m$ and write the resulting map as $\widetilde{w}$.
Over the punctured disc $$
\Sigma = D^2 \setminus (B_1(\frac{\epsilon}{2}) \cup
\cdots B_m(\frac{\epsilon}{2})),$$ the deformed map $\widetilde{w}$ 
does not intersect with the hyperplanes, and it intersects with the 
Clifford torus along the boundaries of the punctured disc.
Recall that the Maslov index is a homotopy invariant. Hence, we have
$\mu(\widetilde{w}|_{D^2}) = \mu(w)$.

Now we will compute the the Maslov index of the map $\widetilde{w}$.
Note that the boundary $\partial \Sigma$ is 
$\partial D^2 \cup ( \cup_{i} \partial B_i(\epsilon/2))$. 

Since the image of the map $\widetilde{w}$ on the
boundaries of the ball $B_i(\frac{\epsilon}{2})$ lies on the Lagrangian submanifold $T^n$,
the map $\widetilde{w}:(\Sigma,\partial \Sigma) \to (\CP^n,T^n)$ satisfies the 
Lagrangian boundary condition. Furthermore, 
since every intersection with the hyperplane
occurs inside the balls $B_i(\epsilon/2)$, $\widetilde{w}|_{\Sigma}$ does not meet the 
hyperplanes. Hence, it can be considered as a map into the cotangent bundle of $T^n$, 
since $\CP^n \setminus (H_0 \cup H_1 \cup \dots \cup H_n)$
can be identified with the cotangent bundle of $T^n$.
 From the Lemma \ref{cot}, 
\begin{equation}\label{indexzero}
\mu(\widetilde{w}|_{\Sigma}) = 0.
\end{equation}
By the Definition in section \ref{maslovindex},  
the Maslov index of the map $\widetilde{w}|_{\Sigma}$ is given by the sum 
of the 
Maslov indices of $\partial \Sigma$ after fixing the trivialization.
 
 Now consider the map $\widetilde{w}:D^2 \to \CP^n$ and we fix a trivialization $\Phi$
 of the pull-back bundle $\widetilde{w}^*T\CP^n$. It gives a trivialization $\Phi_{\Sigma}$
 of the pull-back bundle $(\widetilde{w}|_{\Sigma})^*T\CP^n$ restricted over $\Sigma$.
In this trivialization, it is easy to see that 
 $$\mu(\Phi_{\Sigma},\partial D^2) = \mu(\Phi,\partial D^2) = \mu(\widetilde{w})
 =\mu(w).$$
Since the boundary of the balls $B_i$ are oriented in the opposite way,
and from the explicit description (\ref{map}) of the deformed map on
the ball $B_i$, we have 
$$\mu(\Phi_{\Sigma},\partial B_i) = - 2(\textrm{sum of intersection multiplicities in $B_i$}).$$

From the equation \ref{indexzero}, we have 
 $$\mu(w) - 2(\textrm{sum of intersection multiplicities }) =0.$$
 \end{proof}

\section{Classification and regularity of the holomorphic discs}\label{sec:classify}
With the Maslov index formula in Theorem \ref{indexformula}, we can completely 
classify all holomorphic discs with boundary lying on the Clifford torus.
Here is our classification Theorem.
\begin{theorem}\label{classify}
Let $w :(D^2,\partial D^2) \to (\CP^n,T^n)$ be a holomorphic disc.
Then, homogeneous coordinate functions of the map $w$ can be chosen
so that they are a finite Blaschke products.

i.e. the map $w$ has homogeneous coordinates $[\G_0(z):\cdots:\G_n(z)]$
such that for all $i=0,1,\cdots,n$, there exists $\mu_i\in \ZZ_+$, 
$\alpha_{i,j} \in int(D^2)$ for $j=1,2,\cdots, \mu_i$ and we can write
$$ \G_i(z) = e^{\theta_i}
\prod_{j=1}^{\mu_i}\frac{z-\alpha_{i,j}}{1-\overline{\alpha_{i,j}}z}$$
where $\cap_{i=0}^n \cup_{j=1}^{\mu_i} \{\alpha_{i,j}\} = \phi.$
And the Maslov index of $w$ is $\sum_{i=0}^n \mu_i$.
\end{theorem}
\begin{proof}
 Any map with these coordinate functions are obviously holomorphic, 
But we need to prove that there does not exist any other holomorphic discs.

First, consider the case that the Maslov index of the holomorphic disc is less than 
2n+2. By Theorem \ref{indexformula}, any disc which intersects all 
$(n+1)$ hyperplanes $H_i$'s will have at least $2(n+1)$.
Hence the image should miss at least one hyperplane, say $H_0$. 
Then the map $w$ can be considered as a holomorphic map from $D^2$ to $U_0 \cong \CC^n$ 
with boundary 
in $(S^1)^n \subset \CC^n$. 
Let $\pi_i:\CC^n \to \CC$ be the projection map onto $i$-th coordinate. 
The composition 
$pi_i \circ w : D^2 \to \CC$ maps the boundary $\partial D^2$ to
$S^1 \subset \CC$. But we have a complete classification of such maps.
Namely, they are given by a finite Blaschke products. This proves the
Theorem in the case that the Maslov index of the disc is less than
$2n+2$. 

Now, we consider the case that the Maslov index of the disc is 
bigger or equal to $2n+2$. If such a disc misses at least one
of the hyperplanes $H_i$s, then one can argue similarly as above.
So we assume that the image of the map $w$ intersects with all $n+1$ hyperplanes.
We label every point of the domain of the intersection of the map $w$ with the
fixed hyperplane $H_0$ as $p_1,p_2,\cdots,p_m \in D^2$. And we denote the intersection
multiplicity at the point $p_i$ as $d_i$. 
Let $u:(D^2,\partial D^2) \to (\CP^n,T^n)$ be the map given by 
$$[\prod_{i=1}^m (\frac{1-\overline{p_i}z}{z-p_i})^{d_i} \G_0(z):\G_1(z):\cdots:\G_n(z)]$$
Note that the multiplication preserves the boundary condition, and resulting map is still 
holomorphic.
But the map $u$ no longer intersects with the hyperplane $H_0$, hence 
previous arguments can be applied.
This proves the theorem and the statement about the Maslov index 
follows from Theorem \ref{indexformula}.
\end{proof}
%
Now, we show the regularity of $J_0$ which 
justifies that we may use the standard complex structure to
compute the Floer cohomology.

\begin{theorem}\label{thm:regular}
In the case of the Clifford torus, the standard
complex structure $J_0$ is regular for the holomorphic discs with 
Maslov index less than $2n+2$.
$i.e.\,\, \textrm{Coker}\, D\overline{\partial}_{J_0} = {0}$.
\end{theorem}
\begin{remark}
For the regularity of discs with the Maslov index $\geq 2n+2$,
see \cite{co}
\end{remark}

\begin{proof}
Let $w:(D^2,\partial D^2) \to (\CP^n,T^n)$ be a holomorphic
disc with Maslov index less than $2n+2$. 
Because of Theorem \ref{indexformula}, we consider it as a map
$$w:(D^2,\partial D^2) \to (\CC^n,(S^1)^n).$$
If we linearize at $w$, we obtain a Riemann-Hilbert Problem (see 
\cite{oh:rhp95} or \cite{oh:fclipd293}).
\begin{equation}\label{dbar}
\left\{ \begin{array}{ll}
   \frac{\partial \xi}{\partial \overline{z}} = 0  & \textrm{in}\,\, D^2 \\
   \xi(z) \in T_{w(z)}(S^1)^n & \textrm{for }\,\, z \in \partial D^2
   \end{array} \right.
\end{equation}
where $\xi: D^2 \to \CC^n$ is a smooth map.

Actually, problem \ref{dbar} is completely seperable into $n$ equations
of one variable of the type : for the projection map onto $i$-th coordinate
$\pi_i:\CC^n \to \CC$,
\begin{equation}\label{dbar1}
\left\{ \begin{array}{ll}
   \frac{\partial \eta}{\partial \overline{z}} = 0  & \textrm{in}\,\, D^2 \\
   \xi(z) \in T_{\pi_i(w)}S^1 & \textrm{for }\,\, z \in \partial D^2
   \end{array} \right.
\end{equation}
Now the theorem immediately follows from the study of 1-dimensional 
Riemann-Hilbert problem with this Lagrangian loop:
Oh(\cite{oh:rhp95}) proved the regularity of holomorphic discs with
partial indices $\geq -1$, and in the 1-dimensional problem, partial
index equals the Maslov index which is non-negative in our case.
This finishes the proof.
\end{proof}


For certain elements, for example, the holomorphic disc given by 
$w(z)=[z^2:1:\cdots:1]$, has a nontrivial automorphism if we
do not put any marked point to the moduli space.
$i.e. \,\, z \to e^{\pi i} z$ gives rise to an $\ZZ/2\ZZ$-automorphism group
for the element $((D^2,0),w(z)) \in \CM(\beta_0^2)$.
After putting one or more marked point, the moduli spaces of
holomorphic discs always have trivial automorphisms if it does not
contain sphere bubbles as shown in Lemma 9.1 in \cite{fooo:lift00}.
In our case, when we are considering
moduli space of holomorphic discs with Maslov index $< 2n+2$,
these discs cannot bubble off a sphere since sphere bubbling can only occurs
for discs with Maslov index at least $2n+2$.

\section{Computation of Floer cohomology with 
the standard spin structure}

Since the standard complex structure $J_0$ is regular, we can compute
explicitly the Floer boundary operators using the classification 
Theorem \ref{classify}.

Here we work with the spectral sequence desribed in Theorem \ref{spectral}.
Then the $E_2$ term of the spectral sequence is given by 
$(H^*(L:\QQ) \otimes e^q)^p$. The boundary map on $E_2$ is given by 
$$\delta_2 = \sum_{\beta,\mu(\beta)=2} \delta_{\beta}.$$
Hence, we compute the boundary map $\delta_2$ for 
cohomology generators of $T^n$.

Now we choose the generators of the singular cohomology of
the Clifford torus. We denote by $L_i$ the boundary cycle of the standard disc
$$b_i =[\underbrace{1:\cdots:1}_{i}:z:1:\cdots:1]$$
for $i=0,1,\cdots n$. Then, the cycles $L_i$ for $i=1,2,\cdots n$
generates $H_1(T^n)$,
And we write $L_i \times L_j$ ($i\neq j$) for the cycle given from the 
boundary of 
$$[\overbrace{\underbrace{1:\cdots:1}_{i}:z:\cdots:1}^{j}:z:\cdots:1].$$
We also define  products $L_{i_1} \times \cdots \times L_{i_k}$ similarly.
These products will give all the generators of $H_*(T^n)$.
Recall that we identify these cycles as an element of cohomology by the
relation (\ref{current}).

Now we can compute the boundary operator $\D_{\beta_i}$ for
$\mu(\beta_i)=2$. first,
$$\D_{\beta_i}<pt> = (\CM_2(\beta_i) \,_{ev_1} \times_i <pt>, ev_0)
$$ Here
\begin{eqnarray*}
 [\CM_2(\beta_i)] &=& ([\widetilde{\CM}(\beta)][\partial D^2_0][\partial D^2_1])
 /PSL(2:\RR)\\
&=& (-1)^n([\partial D^2_0][\widetilde{\CM}_(\beta)][\partial D^2_1])/PSL(2:\RR)\\
&=& (-1)^n[\partial D^2_0][T^n] \;\;\;\textrm{from Theorem \ref{tn}}
\end{eqnarray*}
where $[\;]$ means oriented basis of the tangent space at any 
element $((D^2,z_0,z_1),w) \in \CM_2(\beta_i)$.
Here $\partial D^2_i$ denotes $i$-th marked point. So by the
definition of the orientation of the fibre product in section \ref{sec:ori},
$$\D_{\beta_i}<pt> = (\CM_2(\beta_i) \,_{ev_1} \times_i <pt>, ev_0)
= (-1)^n L_i$$
where we obtained $L_i$ as we evaluate the marked point $z_0$ along the
boundary of the disc $D_0^2$.
And similarly for ($i\neq j$)
\begin{equation}\label{d2}
\D_{\beta_j} L_i = (\CM_2(\beta_j)\,_{ev_1}
 \times_i L_i, ev_0)
  =(-1)^n ( L_j \times L_i) 
\end{equation}
 For $i=j$, one can easily see that $\D_{\beta_i} L_i = 0$. Hence
 we will use the equation (\ref{d2}) even for $i=j$ with the convention
 $L_i \times L_i =0$.

If we take a sum of $\D_{\beta_i}$ for all $i=0,1,\cdots
n$, then
\begin{eqnarray*}
\D_2 <pt> &=& (-1)^n (L_0 +
\cdots + L_n) \\
&\equiv&  (-1)^n ((- L_1 - \cdots - L_n ) + L_1 +  \cdots + L_n) \;=\;0 \;\;\; 
\textrm{in}\; H^*(T^n:\QQ)
\end{eqnarray*}
For the higher dimensional generators, we can proceed similarly.
For any generator of the  singular cohomology of $T^n$ represented as
$L_{i_i}\times L_{i_2} \times \cdots \times L_{i_k}$, we compute the
boundary of it as
$$ \D_2 (L_{i_i}\times L_{i_2} \times \cdots \times L_{i_k}) =$$
$$(-1)^n(L_0+L_1+\cdots+L_n)\times (L_{i_i}\times L_{i_2} \times \cdots \times L_{i_k})
\equiv 0\;\;\; 
\textrm{in}\; H^*(T^n:\QQ)$$
We can see that under the standard spin-structure, the
boundary operator $\delta_2 \equiv 0$ for all cohomology generators of $T^n$.

For the boundary operators $\delta_k$ for $k \geq 4$, we have the 
following Proposition, which follows from the description of the
moduli space of holomorphic discs.
\begin{prop}\label{zero4}
For the Bott-Morse Floer cohomology of the Clifford torus, the
boundary operator $\delta_\beta \equiv 0 $ for $\mu(\beta) \geq 4$
\end{prop}
\begin{proof}
We will show that when $\mu(\beta) \geq 4$, the dimension of the image
under the evaluation map $ev_0$ of the moduli space
$(\CM_2(\beta)\,_{ev_1} \times_f P)$ is always less than the 
dimension of moduli space itself. Hence it will prove that
$\delta_{\beta} \equiv 0$ as we consider them as currents.

Consider any homotopy classes $\beta \in \pi_2(\CP^n,T^n)$ with 
$\mu(\beta) =4$. The dimension of the moduli space 
$(\CM_2(\beta)\,_{ev_1} \times_f P)$ is $dim(P) +3$. 
But we claim that, for any point $<pt> \in P$, 
$$dim(ev_0(\CM_2(\beta)\,_{ev_1} \times_f <pt> )) \leq 2.$$
The claim  easily follows from the classification theorem:
We argue by example. Consider the homotopy class $(\beta_0+\beta_1)$,
which is the homotopy class of the map 
$$[\frac{z-\alpha_1}{1-\overline{\alpha_1}z}:
\frac{z-\alpha_2}{1-\overline{\alpha_2}z} :1:\cdots:1]$$ 
where $\alpha_i \in D^2$ for $i=1,2$.
$$ev_0((\CM_2(\beta_0 +\beta_1)\,_{ev_1} \times_f <pt> ) \subset
[p_0 e^{i\theta_1}: p_1e^{i\theta_2} : p_2 :\cdots : p_n]$$
where $[p_0:\cdots:p_n]$ represents point $f(p)$ in $T^n$
and $0\leq \theta_i \leq 2\pi$ for $i=1,2$.
Hence 
the dimension of the chain $(\CM_2(\beta) \times P,ev_0)$ is at most
$dim(P) +2$. Hence, it gives zero as a current. $i.e. \,\,\delta_{\beta} \equiv 0$.

The above argument can be easily generalized for 
homotopy classes with higher Maslov indices.
\end{proof}

Since all quantum boundary operators are zero, 
the spectral sequence degenerates at $E_2$. Hence we have 

%
\begin{theorem}\label{stan}
For the standard spin structure, we have an isomorphism of 
$\Lambda_{nov}$-modules with $\ZZ$-grading. 
$$ HF^*(T^n,T^n;\Lambda_{nov}) \cong H^*(T^n) \otimes \Lambda_{nov}.$$
\end{theorem}
\begin{remark}
This isomorphism does not preserve the product structure. It is only
a module isomorphism.
\end{remark}

\section{The Floer cohomology with non-standard spin structures.}

In this section, we compute the Floer cohomology with non-standard spin 
structures. It will be done by describing the change of sign in the 
Floer coboundary operator according to the change of spin structure.
As we have described in Theorem \ref{sgnchange},
the change of  spin-structure of $T^n$ results in 
the change of orientations of certain moduli spaces $\CM(\beta)$
for $\beta \in \pi_2(M,L)$.
There exist $|H^1(T^n;\ZZ/2\ZZ)| = 2^n$ spin structures of
the Clifford torus. We denote elements of $\ZZ/2\ZZ$ as 0 and
1. We label the standard spin structure as $(0,0,\cdots,0) \in
(\ZZ/2\ZZ)^n$. Let $I$ be a subset of $\{1,2,\dots,n\}$. Then
Spin$_I$ will denote the spin structure corresponding to
$(a_1,a_2,\dots,a_n) \in (\ZZ/2\ZZ)^n$ where $a_i = 1$ for $i \in
I$ and $a_i =0$ for $i \notin I$.
Let Spin$_0$ denote
the standard spin structure.

But we will give a different labeling of the spin structures as
follows. Let $\varepsilon_i \in \{-1,+1\}$ for $i \in
\{0,1,\cdots,n\}$
\begin{definition}  
Consider $(\varepsilon_0,\cdots,\varepsilon_n) \in
\{-1,+1\}^{n+1}$ which satisfies
$$ \varepsilon_0 \cdot \varepsilon_1 \cdots \varepsilon_n =1$$
And let $$I := \{ i  | \varepsilon_i =-1, i\neq 0 \}$$ Then
 we will denote the spin structure $Spin_I$ by 
 $(\varepsilon_0,\cdots,\varepsilon_n)$.
\end{definition}
This labeling is more convenient because $\varepsilon_i$
will be the orientation change of the moduli space
$\widetilde{\CM}_2(\beta_i)$, for each $i =0,1,\cdots n$ when we
change the
standard spin structure to the spin structure Spin$_I$:
\begin{equation}\label{change}
[\CM(\beta_i)]_{Spin_I} = \varepsilon_i [\CM(\beta_i)]_{Spin_0}
\end{equation}
The reason is that for $i\neq 0$, if $i\in I$, the spin structure
Spin$_I$ will twist the trivialization of tangent bundle of $T^n$
along the $L_i$ from the standard trivialization, which will
change the orientation of the moduli space
$\widetilde{\CM}_2(\beta_i)$. So $\varepsilon_i$ is exactly the
sign change of the moduli space $\widetilde{\CM}_2(\beta_i)$, or
the sign change of the boundary operator $\delta_{\beta_i}$. When
$i =0$, orientation change of $\widetilde{\CM}_2(\beta_0)$ will
depend on the product $ \varepsilon_1 \cdots \varepsilon_n=
\varepsilon_0$ since the boundary of $\beta_0$ disc has homology
class $L_0 \cong -L_1-L_2-\cdots L_n$ in $\pi_1(T^n)$.

Similar sign changes occur for homotopy classes 
with higher Maslov indices according to their boundary elements.
(according to the map $(\pi_2(M,L) \to \pi_1 (L))$
But because of Proposition \ref{zero4}, it will be irrelevant to
the Floer cohomology.

Now we calculate the Floer cohomology of $T^n$ with these spin
structures. We fix our spin structure by
$(\varepsilon_0,\cdots,\varepsilon_n)$ with $ \varepsilon_0 \cdot
\varepsilon_1 \cdots \varepsilon_n =1$ or Spin$_I$.
We consider the homotopy classes with Maslov index 2. 
Recall that they are indexed as $\beta_i$ for $i =0,\dots n$.
From the sign change rule (\ref{change}), we have 
\begin{eqnarray}\label{delspin}
\delta_2 <pt> &=&\sum_{i=0}^n \delta_{\beta_i} <pt>\; = \;(-1)^n \sum_{i=0}^n \varepsilon_i L_i 
\nonumber \\
&=& (-1)^n \sum_{i=1}^n (\varepsilon_i - \varepsilon_0) L_i 
\;\;\; 
\textrm{in}\; H^*(T^n:\QQ)
\end{eqnarray}
Last equality can be obtained by writing $L_0$ 
as $(-L_1-L_2- \cdots-L_n)$. Hence,
$\delta_\beta<pt> =0$ if and only if $\varepsilon_i=\varepsilon_0$ for all
$i$. 
There exist at most two such spin structures. 
The case $\varepsilon_0=\varepsilon_1=\cdots =\varepsilon_n =
1$ is the standard spin structure case, and the case
$\varepsilon_0=\varepsilon_1=\cdots =\varepsilon_n = -1$ can occur
only when $n$ is an odd integer because it has to satisfy $
\varepsilon_0 \cdot \varepsilon_1 \cdots \varepsilon_n =1$. Except
these two possibilities, $\delta_2 <pt> \neq 0$. 
Therefore, $<pt>$ is no longer in the kernel of $\delta$.

Later, we will consider Floer cohomology with flat line bundle $\CL$ on $T^n$.
Let $h_i$ be the holonomy of $\CL$ along the cycle $L_i$.
Then Floer cohomology with spin structure $(\varepsilon_0,\cdots,\varepsilon_n)$
is same as that with flat line bundle $\CL$ whose holonomy is given as
\begin{equation}
h_i = \varepsilon_i.
\end{equation}
(Compare (\ref{delspin}) with (\ref{delhol}) for example).

\begin{remark}
These various spin-structures can be realized as being
given a flat complex line bundle with holonomy $e^{\pi i}$ for the
corresponding generators. More precisely, $H^1(L,\ZZ/2\ZZ)$ which
characterizes the spin structures, also gives the flat
real line bundles over $L$. We get the corresponding flat complex
line bundle by tensoring $\CC$ to this real line bundle.
\end{remark}
Then, the Theorem \ref{dbranefl} can be interpreted in terms of spin structures 
as follows:
\begin{theorem}
For $n$ even, with any non-standard spin structure, 
Floer cohomology $H^*(T^n,T^n;\Lambda_{nov})$ 
vanishes. 

For $n$ odd, let $(0,\dots 0) \in (\ZZ/2)^n$ be the standard spin structure.
Then Floer cohomology for other spin structures vanishes
except the spin structure $ (1,1,\dots,1) \in (\ZZ/2\ZZ)^n$, 
in which case
$$ HF^*(T^n,T^n;\Lambda_{nov}) \cong H^*(T^n) \otimes \Lambda_{nov}.$$
\end{theorem}
\begin{proof}
When $n$ is even, non-standard spin structures does not give 
specified holonomies whose Floer cohomology is non-vanishing.
But when $n$ is odd, for $k = (n+1)/2$, 
\begin{equation}
e^{\frac{2\pi k}{n+1}} = e^{\pi i}.
\end{equation}
Then the theorem follows from Theorem \ref{dbranefl}
\end{proof}

\begin{remark}
We may define Floer cohomology with Novikov ring with $\ZZ$ coefficient
since $T^n$ is monotone. In this coefficient ring, with non-standard spin structure, it gives
a non-vanishing Floer cohomology becasue it will have a torsion
element. (We can not divide by 2 in (\ref{divide}) for example ).
\end{remark}

\section{D-branes and Floer cohomology}
The following definitions are from \cite{fukaya:fhms99}.
\begin{definition}{\cite{fukaya:fhms99}}
Let $(M,\omega)$ be a symplectic manifold. A pair
$(L,\mathcal{L})$ of Lagrangian submanifold $L$ of $M$ and a flat
complex line bundle $\mathcal{L}$ on L is a \textbf{brane} (in a
classical sense) of A-model compactified by $(M,\omega)$.
\end{definition}
Floer cohomology of the above pairs
$(L_0,\mathcal{L}_0),(L_1,\mathcal{L}_1)$ was proposed by 
Konsevich~\cite{konsevich:hams95}. One
can define a Bott-Morse Floer cohomology of the pair
$(L_0,\mathcal{L}_0)$ by modifying the boundary operator as
follows. And we 
use the Novikove ring with $\CC$-coefficient $\Lambda_{\CC,nov}$
instead of $\QQ$-coefficient (see Definition \ref{novikov}).
\begin{definition} 
We define Bott-Morse D-brane Floer cohomology 
of the cochain complex $C^*(L,\Lambda_{\CC,nov})$
by defining the coboundary map
as 
\begin{equation}
\begin{cases}
\D_{\beta}([P,f])=(\CM_2(\beta) \,_{ev_1}\times_f P, ev_0)\cdot
 (hol_{\partial \beta} \mathcal{L}) \,\,\,\textrm{for }\, \beta \neq 0,\\
\D_0([P,f])=(-1)^n[\partial P,f]
\end{cases}
\end{equation}
where $hol_{\partial \beta}\mathcal{L}$ is the holonomy of the
flat line bundle along the closed curve $\partial \beta$.
And define the coboundary map $\delta$ as 
$$\delta([P,f])=\sum_{\beta \in \pi_2(M,L)}
\delta_\beta([P,f]) \otimes e^{\frac{\mu(\beta)}{2}}.$$
\end{definition}
Then, $\delta \circ \delta =0$ follows from Theorem 
\ref{squarezero}.
Note that $\delta_0$ does not change from definition
of Bott-Morse Floer cohomology.
It corresponds to the fact that thin-trajectories in Floer
cohomology between $L$ and $\phi(L)$
 will not have any holonomy factor if we consider the flat line
 bundle induced on $\phi(L)$ from $\mathcal{L}$ .

Now, let $h_j \in S^1$ denote the holonomy of $\mathcal{L}$
along the generators $L_j$. Then we can see that
$$\D_{\beta_j}<pt>= ( (-1)^n L_j) \cdot h_j$$
$$\D_{\beta_0}<pt>= ( (-1)^n L_0) \cdot
h_1^{-1}\cdot h_2^{-2}\cdots h_n^{-1} $$ 
Consider the case that $$h_j =
e^{\frac{2\pi k i}{n+1}} \,\, \textrm{for all} \, j = 1,2,\cdots,n
, \,\, \textrm{for a fixed} \,\,k \in \ZZ $$ 
Then $$h_0 =
h_1^{-1}\cdot h_2^{-2}\cdots h_n^{-1} = e^{-\frac{2\pi n k
i}{n+1}} = e^{\frac{2\pi k i}{n+1}} = h_j$$ 
Therefore, in this case, boundary
operators $\delta_{\beta_j}(<pt>)$ are multiplied by the factor $e^{\frac{2\pi k
i}{n+1}}$ for all $j$, so under the standard spin-structure,
the boundary operator $\delta_2$ is still zero. Hence we get the same
Floer cohomology as in Theorem \ref{stan}.
\begin{theorem}\label{dbranefl}
Let the following
$$(1,\cdots,1),(\alpha,\cdots,\alpha),\cdots,(\alpha^{n},\cdots,\alpha^{n})$$
for $\alpha = e^{\frac{2\pi i}{n+1}}$ represent the holonomies of
the flat line bundles along the generators of the Clifford torus
$T^n$ for $n \geq 1$. Under the standard spin-structure, the above A-branes are
the only ones which give non-trivial Floer cohomology, which are
isomorphic to the singular cohomology as in Theorem \ref{stan}
\end{theorem}
\begin{remark}
The above Theorem confirms the prediction of Hori~\cite{hori:lmsd00}.
See Kapustin-Li~\cite{kapustin:dlgmag02}
for explicit statement on the product structure for the case $n=2$.
\end{remark}
\begin{proof}
To finish the proof of the above Theorem, we need to prove that for
flat bundles with other holonomies, D-brane Floer cohomology vanishes. It
will be similar to the calculation of Floer cohomology with a non-standard spin
structure.
Again, we work with the spectral sequence and we start with $E_2$ as before.
Recall that $h_j$ denotes the holonomy along the generator $L_j$.
Define 
\begin{equation}
\begin{cases}
 S:= \{j\in \{1,2,\cdots,n\}| h_j = h_0 \} \\
 S^c:= \{j\in \{1,2,\cdots,n\}| h_j \neq h_0\}
\end{cases}
\end{equation}

First, the case when $S^c$ is empty is exactly
the case when $h_j = e^{\frac{2\pi k i}{n+1}}$ for all
$j=0,1,\cdots,n$ for some $k$. Then $\delta$ is always zero, hence
Floer cohomology has the same generator as the singular cohomology
of
$T^n$ as in the Theorem.

So let us assume that $S^c$ is not empty. Then
\begin{eqnarray}\label{delhol}
\delta_2 <pt> &=&\sum_{j=0}^n \delta_{\beta_j} <pt> \;= \;
(-1)^n \sum_{j=0}^n h_j L_j   \nonumber \\
&=& (-1)^n \sum_{j=1}^n (h_j - h_0) L_j  
\;=\; (-1)^n \sum_{j\in S^c} (h_j - h_0) L_j  
\end{eqnarray}

So $<pt>$ is not in the kernel of $\D_2$. For 1 dimensional cycles,
\begin{eqnarray*}
\delta_2 &(&\sum_{j=1}^n a_j L_j) =
\sum_{i=0}^n\sum_{j=1}^n a_j ( \delta_{\beta_i} L_j)\\
&=& (-1)^n \sum_{i=0}^n \sum_{j=1}^n a_j(h_i L_i \times L_j)\\
&=& (-1)^n \sum_{i=1}^n \sum_{j=1}^n a_j(h_i - h_0) L_i \times L_j \\
&=& (-1)^n \sum_{1\leq i<j \leq n}L_i \times L_j (a_j(h_i - h_0) -
a_i(h_j - h_0))
\end{eqnarray*}
So to be in the kernel,
\begin{equation}\label{diff2}
(a_j(h_i - h_0) -
a_i(h_j - h_0))=0
\end{equation}
for all $i,j$ with $i < j$.
Now, for $i,j \in S^c$, the equation \ref{diff2} becomes
$$\frac{a_j}{(h_j-h_0)} = \frac{a_i}{(h_i-h_0)} $$
So for all  $j\in S^c$ , $\frac{a_j}{(h_j-h_0)}$ are same
and  we denote it by $a$.
For $i \in S$, $j \in S^c$ with $i<j$, the equation 
\ref{diff2} becomes $a_i=0$.
Similarly, for $i \in S$, $j \in S^c$ with $i>j$, the equation
\ref{diff2} implies $a_i=0$.

Now the elements in the kernel can be written as
$$ \sum_{i \in S^c} a_i L_i = \sum_{i \in S^c} a (h_i - h_0) L_i$$
In fact this element is in the image of $\delta_2$:
One can check that 
$$\delta_2 \;((-1)^n a <pt> ) 
=  \sum_{i \in S^c} a (h_i - h_0) \,L_i$$
Similarly, for higher dimensional cycles, we show that the elements in the 
kernel of $\delta_2$ lies in the image of $\delta_2$.

First, we will set up our notation for the indices
\begin{definition}
By $G,I,J$ we will denote subsets of $\{1,2,\cdots,n\}$ with number of 
elements $|G|=k-1,|I|=k,|J|=k+1$. We will also denote its 
elements as $G=\{g_1,\cdots,g_{k-1}\}$ with 
$g_1<g_2<\cdots<g_{k-1}$. And we denote $G_{\widehat{s}} =
G \setminus \{g_s\}$. This notation will be applied to any index set.
\end{definition}

We need the following elementary lemma, which states that 
the $k$-th simplicial cohomology of the standard $(n-1)$-simplex is zero.
Let $R$ be a coefficient ring which will be either $\QQ$ or $\CC$.
\begin{lemma}\label{cech}
We fix $1 \leq k \leq n$.
Suppose we are given numbers $A_I \in R$ for every subset $I$ with $|I|=k$.
And suppose those numbers satisfy the following equation: 
For any J with $|J|=k+1$, 
$$\sum_{s=1}^{k+1} (-1)^{s-1}A_{J_{\widehat{s}}} =0$$

Then, there exists $B_G \in R$ for all $G$ with $|G|=k-1$ so that 
each $A_I$ can be written as 
$$A_I = \sum_{s=1}^{k} (-1)^{s-1}B_{I_{\widehat{s}}}.$$
When $k=1$, there exist a number $B_G$ for the empty set $G$ with $|G|=0$
so that 
$$A_I = \sum_{s=1}^{k}(-1)^{s-1}B_G.$$
\end{lemma}
\begin{proof}
There is an obvious correspondence between the index set, say 
$I \subset \{1,2,\cdots,n\}$, and 
the simplicial chain, say $C_I$, of the standard $(n-1)$-simplex.
Then, consider $A_*$ as a simplicial $k$ cochain which assigns the number $A_I$ to
the chain $C_I$.
Then the hypothesis is equivalent to the fact that $\delta A_* =0$.
Hence, there exists $k-1$ dimensional cochain $B_*$ with 
$$\delta B_* = A_*$$.
\end{proof}

Now we denote an arbitrary element of $k$ dimensional cycles as 
$$\sum_{I,|I|=k} A_I L_I$$
The boundary of this element is
\begin{eqnarray*}
\delta_2 (\sum_{I,|I|=k} A_I L_I) &=& \sum_{I} A_I (\delta L_I) \\
&=& \sum_{I} A_I (h_0 L_0 + \cdots + h_n L_n ) \times L_I \\
&=& \sum_I \sum_{s=1}^n (h_s - h_0)A_I L_s \times L_I \\
&=& \sum_{J,|J|=k+1}\sum_{s=1}^{k+1} A_{J_{\widehat{s}}}(-1)^{s-1} (h_{j_s} - h_0) L_J
\end{eqnarray*}

Hence, the element $\delta_2 (\sum_{I,|I|=k} A_I L_I)$ is in the kernel 
if for all $J$ with $|J|=k+1$,
\begin{equation}\label{kernel}
\sum_{s=1}^{k+1} A_{J_{\widehat{s}}}(-1)^{s-1} (h_{j_s} - h_0) =0
\end{equation}

Now we show that any element in the kernel lies in the image of
the boundary map.
First, we consider the case that 
the set $ S= \{i\in \{1,2,\cdots,n\}| h_i = h_0\}$ is empty.
i.e. $h_0 - h_i \neq 0$ for all $i =1,\cdots,n$.
For all $I, |I|=k$, we set 
\begin{equation}\label{divide}
B_I = \frac{A_I}{\prod_{i\in I}(h_i-h_0)}.
\end{equation}
Then, the equation (\ref{kernel}) is nothing but
$$\sum_{s=1}^{k+1} (-1)^{s-1}B_{J_{\widehat{s}}} =0$$

By the Lemma \ref{cech}, there exists $C_G$ for all $G\subset \{1,2,\cdots,n\}$ with
$|G| = k-1$ such that
$$ B_I =  \sum_{s=1}^{k} (-1)^{s-1}C_{I_{\widehat{s}}}$$

Then we take the boundary of $\sum_G (\prod_{i\in G}(h_i-h_0))C_G L_G$, where
the sum is taken over all $G\subset \{1,2,\cdots,n\}$ with $|G| =k-1$.
\begin{eqnarray*}
\delta_2 (\sum_G (\prod_{i\in G}(h_i-h_0))C_G L_G) &=& \sum_G\sum_{s=0}^n 
C_G(\prod_{i\in G}(h_i-h_0))(h_s-h_0)L_s L_G \\
&=& \sum_I (\prod_{i\in I}(h_i-h_0))B_I L_I \\
&=& \sum_I A_I L_I
\end{eqnarray*}
This shows that codimension k kernels are coboundaries.

Now when $S$ is not empty, 
we carry out the same argument, but it will be more complicated to prove it.
Without loss of generality, we may set 
$$S=\{r+1,r+2,\cdots,n\}, \;\;\;  S^c = \{1,2,\cdots,r\}$$
From now on, we write $$ A_I = A_{I\cap S^c}^{I \cap S}.$$
For example, we write $A_{1,2,\cdots,r+2}$ as $A_{1,2,\cdots,r}^{r+1,r+2}$.
This is to distinguish elements in $S$ and $S^c$.
Now the codimension k element in the kernel  can be written as 
$(\sum_{I,|I|=k} A_I L_I)$, where $A_I$ satisfies the equation
 (\ref{kernel}).
But if $j_s \in S$, we have $h_{j_s} -h_0 =0$. 
Hence we may write, for each fixed $J$ with $|J|=k+1$,
\begin{equation}\label{k2}
\sum_{j_s \in S^c} A_{J_{\widehat{s}}\cap S^c}^{J\cap S}
(-1)^{s-1} (h_{j_s} - h_0) =0 
\end{equation}

Now, for $P \subset S^c, T \subset S$, we let 
$$B_P^T = \frac{A_P^T}{\prod_{p\in P}(h_p-h_0)}$$
Then, the equation (\ref{k2}) is equivalent to 
$$\sum_{s=1}^{|J\cap S^c|} (-1)^{s-1}B_{(J\cap S^c)_{\widehat{s}}}^{J\cap S} =0$$
We collect all such equations with respect to the same index set $J\cap S$.
By the Lemma \ref{cech}, there exists $C_Q^{J\cap S}$ 
for all $Q\subset S^c$ with
$|Q| = |J\cap S^c|-2$ such that
for any $P \subset S^c$ with $|P|=|J\cap S^c|-1$ 
$$ B_P^{J \cap S}=  \sum_{s=1}^{|P|} (-1)^{s-1}C_{P_{\widehat{s}}}^{J \cap S}$$
Then we may rewrite the above as 
$$ B_{I\cap S^c}^{I \cap S} = \sum_{s=1}^{|I\cap S^c|}(-1)^{s-1}
C_{(I\cap S^c)_{\widehat{s}}}^{I \cap S} $$
Now we will show that $\sum_I A_I L_I$ is exactly in the image of the 
following element.
In the following, we take the sum over any subset $T \subset S$, and over
any subset $Q \subset S^c$ with $|T| + |Q| = |I| -1$ and 
for any subset $P \subset S^c$ with $|T| + |P| = |I|$. 

\begin{eqnarray*}
\delta_2(\sum_T\;\sum_Q(\prod_{q\in Q}(h_q-h_0))
\,C_Q^T \,L_Q \,L_T) &=&
\sum_T\;\sum_Q (\prod_{q\in Q}(h_q-h_0)) 
\,C_Q^T\,(\sum_{s=1}^r(h_s-h_0)\,L_s \,L_Q L_T)) \\
&=& \sum_T\;\sum_P (\prod_{p\in P}(h_p-h_0))\,
B_P^T \,L_P \,L_T \\
&=&  \sum_T\;\sum_P A_P^T\, L_P\, L_T \\
&=& \sum_T \;\sum_P A_{P\cup T}\, L_{P \cup T} \\
&=& \sum_I A_I\, L_I 
\end{eqnarray*}
This finishes the proof. 
\end{proof}

\bibliographystyle{amsalpha}

\end{document}